\documentclass[journal,twoside,web]{ieeecolor}
\usepackage{generic}
\usepackage{cite}
\usepackage{amsmath,amssymb,amsfonts}
\usepackage{algorithmic}
\usepackage{graphicx}
\usepackage{textcomp}
\def\BibTeX{{\rm B\kern-.05em{\sc i\kern-.025em b}\kern-.08em
    T\kern-.1667em\lower.7ex\hbox{E}\kern-.125emX}}
\usepackage{cite}
\usepackage{amsmath,amssymb,amsfonts}
\usepackage{algorithmic}
\usepackage{textcomp,float}
\usepackage{graphicx,subfigure}
\usepackage{epsfig} 
\usepackage{amsmath} 
\usepackage{amssymb}  
\usepackage{amsfonts}
\usepackage{euscript}
\usepackage{color}
\usepackage{ccaption}
\usepackage{listings}
\usepackage{booktabs} 
\usepackage{longtable}
\usepackage{multirow}
\usepackage{graphicx}
\usepackage{epstopdf}
\usepackage{bm}
\usepackage{makecell}
\usepackage{subfig}
\usepackage{enumerate}

\newtheorem{theorem}{\bf \text{Theorem}}

\newtheorem{assumption}{\bf \text{Assumption}}

\newtheorem{lemma}{\bf \text{Lemma}}
\newtheorem{remark}{\bf \text{Remark}}

\newcommand{\Z}{{\mathbb Z}}
\newcommand{\E}{{\mathbb E}}
\newcommand{\COV}{{\mathbb C\mathbb O\mathbb V}}

\newcommand{\R}{{\mathbb R}}
\newcommand{\N}{{\mathbb N}}

\newcommand{\TR}{\text{R}}

\makeatletter

\newcommand{\Rmnum}[1]{\expandafter \@slowromancap\romannumeral #1@}
\makeatother

\DeclareMathOperator*{\rank}{rank}

\DeclareMathOperator*{\cond}{cond}
\DeclareMathOperator*{\Tr}{Tr}

\DeclareMathOperator*{\LS}{LS}

\DeclareMathOperator*\argmin{arg\,min}

\def\qed{\hfill \rule{4pt}{7pt}}
\def\pf{\textit{Proof. }}

\begin{document}
\title{A Tutorial on Asymptotic Properties of Regularized Least Squares Estimator for Finite Impulse Response Model}
\author{Yue Ju$^{1}$, Tianshi Chen$^{1,*}$, Biqiang Mu$^{2}$ and Lennart Ljung$^{3}$
	\thanks{*This work was supported in part by the Thousand Youth Talents Plan funded by the central government of China, the general project funded by NSFC under contract No. 61773329, the Shenzhen Science and Technology Innovation Council under contract No. Ji-20170189 (JCY20170411102101881), the Robotic Discipline Development Fund (2016-1418) from Shenzhen Government,  the President’s grant under contract No. PF. 01.000249, the Start-up grant under contract No. 2014.0003.23 funded by CUHKSZ, and by the Swedish Research Council, contract 2019-04956 and by Vinnova’s center LINKSIC.}
	\thanks{$^{1}$Yue Ju and Tianshi Chen* (corresponding author) are with the School of Data Science and Shenzhen Research Institute of Big Data, The Chinese University of Hong Kong, Shenzhen, 518172, China, {\tt\small yueju@link.cuhk.edu.cn, tschen@cuhk.edu.cn}}%
	\thanks{$^{2}$Biqiang Mu is with Key Laboratory of Systems and Control, Institute of Systems Science, Academy of Mathematics and System Science, Chinese Academy of Sciences, Beijing 100190, China {\tt\small bqmu@amss.ac.cn}}
	\thanks{$^{3}$Lennart Ljung is with the Division of Automatic Control, Department of Electrical Engineering, Link$\ddot{\text{o}}$ping University, Link$\ddot{\text{o}}$ping SE-58183, Sweden {\tt\small ljung@isy.liu.se}}
}

\maketitle

\begin{abstract}
In this paper, we give a tutorial on asymptotic properties of the Least Square (LS) and Regularized Least Squares (RLS) estimators for the finite impulse response model with filtered white noise inputs. We provide three perspectives: the almost sure convergence, the convergence in distribution and the boundedness in probability. On one hand, these properties deepen our understanding of the LS and RLS estimators. On the other hand, we can use them as tools to investigate asymptotic properties of other estimators, such as various hyper-parameter estimators.
\end{abstract}




\section{Least Squares Estimator for Finite Impulse Response Model}\label{sec:ls estimator for fir model}

We focus on the $n$th-order finite impulse response (FIR) model as follows,
\begin{align}\label{eq:liner regression model at time t}
y(t)=\sum_{i=1}^{n}g_{i}u(t-i)+v(t),\ t=1,\cdots,N
\end{align}
where $n$ is the order of the FIR model, $N$ denotes the sample size, $t$ denotes the time index,  $u(t)\in\R$, $y(t)\in\R$, and $v(t)\in\R$ are the input, output and disturbance at time $t$, respectively, and $g_{1},\cdots,g_{n}\in\R$ are FIR model parameters to be estimated.

\begin{assumption}\label{asp:input}
	The input $u(t)$ with $t=1-n,\cdots,N-1$ is the filtered white noise with the stable filter $H(q)$, i.e.
	\begin{subequations}\label{eq:asp of ut}
		\begin{align}
		H(q)=&\sum_{k=0}^{\infty}h(k)q^{-k}\ \text{with}\ \sum_{k=0}^{\infty}|h(k)|<\infty,\\
		u(t)=&H(q)e(t)=\sum_{k=0}^{\infty}h(k)e(t-k),
		\end{align}
	\end{subequations}
	where $q^{-1}$ represents the backward shift operator: $q^{-1}u(t)=u(t-1)$, and $e(t)$ is independent and identically distributed (i.i.d.) with zero mean, variance $\sigma_{e}^{2}>0$, bounded moments of order $4+\delta$ for some $\delta>0$, and $\E[e^4(t)]=c\sigma^4_{e}$ with $c\in\R$ being a constant. Moreover, let $\Sigma\in\R^{n\times n}$ be 
	\begin{align}\label{eq:def of Sigma}
	\Sigma\triangleq\COV\left(\left[\begin{array}{cccc}u(0)&u(1)&\cdots&u(n-1)\end{array}\right]^{T}\right),
	\end{align}
	where $\COV(\cdot)$ denotes the covariance matrix of a random vector. In addition, $\Sigma$ is assumed to be positive definite, i.e. $\Sigma\succ 0$.
\end{assumption}

\begin{assumption}\label{asp:noise}
	The measurement noise $v(t)$ with  is i.i.d. with zero mean, variance $\sigma^2>0$ and bounded moments of order $4+\delta$ for some $\delta>0$.
\end{assumption}

\begin{assumption}\label{asp:independence between input and noise}
	$\{u(t)\}_{t=1-n}^{N-1}$ and $\{v(t)\}_{t=1}^{N}$ are mutually independent, which means that for $i=1-n,\cdots,N-1$ and $j=1,\cdots,N$, $u(i)$ and $v(j)$ are independent.
\end{assumption}

By Assumption \ref{asp:input}, $u(t)$ is a stationary stochastic process with
\begin{subequations}\label{eq:asp of vt}
	\begin{align}
	\label{eq:zero mean of ut}
	\E[u(t)]&=0,\\
	\label{eq:definition of R_u_tau}
	\E[u(t)u(t+\tau)]&\triangleq R_{u}(\tau)=\sigma_{e}^2\sum_{k=0}^{\infty}h(k)h(k+\tau),
	\end{align}
\end{subequations}
where $\E(\cdot)$ denotes the mathematical expectation, $\tau\geq0$, and $R_{u}(\tau)=R_{u}(-\tau)$, and moreover, the $(i,j)$th element of $\Sigma$ is $R_{u}(|i-j|)$. 

The model \eqref{eq:liner regression model at time t} can also be rewritten in matrix form as
\begin{align}\label{eq:vector-matrix form linear regression model}
Y=\Phi\theta+V,
\end{align}
where
\begin{subequations}
\begin{align}
Y=&\left[\begin{array}{cccc}y(1)& y(2)& \cdots & y(N)\end{array}\right]^{T},\\
\Phi=&\left[\begin{array}{cccc}\phi(1)& \phi(2) & \cdots & \phi(N)\end{array}\right]^{T},\\
\theta=&\left[\begin{array}{cccc}g_{1}& g_{2} & \cdots & g_{n}\end{array}\right]^{T},\\
V=&\left[\begin{array}{cccc}v(1)& v(2)& \cdots & v(N) \end{array}\right]^{T}
\end{align}
\end{subequations}
with $\phi(t)=\left[\begin{array}{cccc}u(t-1)& u(t-2) & \cdots&u(t-n)\end{array} \right]^{T}$ and $u(t)=0$ for $t< 0$.

Assume that $\Phi\in\R^{N\times n}$ with $N>n$ is full column rank, i.e. $\rank(\Phi)=n$. To estimate the unknown $\theta$ based on data $\{y(t),\phi(t)\}_{t=1}^{N}$, one classic estimation method is the Least Squares (LS):
\begin{subequations}\label{eq:LS estimator}
	\begin{align}\label{eq:LS form1}
	\hat{\theta}^{\LS}=&\argmin_{\theta\in\R^{n}}\|Y-\Phi\theta\|_{2}^{2}\\
    \label{eq:LS form2}
	=&(\Phi^{T}\Phi)^{-1}\Phi^{T}Y,
	\end{align}
\end{subequations}
where $\|\cdot\|_{2}$ denotes the Frobenius norm.
The LS estimator has some interesting asymptotic properties, which will be discussed in the following context.


\section{Asymptotic properties of LS estimator for FIR model}

Let the true FIR parameter be $\theta_{0}\in\R^{n}$. Based on \eqref{eq:LS form2}, the LS estimator can be rewritten as
\begin{align}
\hat{\theta}^{\LS}
=&(\Phi^{T}\Phi)^{-1}\Phi^{T}(\Phi\theta_{0}+V)\nonumber\\
\label{eq:LS estimator form3}
=&\theta_{0}+N(\Phi^{T}\Phi)^{-1}\frac{\Phi^{T}V}{N},
\end{align}
 which consists of three building blocks: $\theta_{0}$, $N(\Phi^{T}\Phi)^{-1}$ and ${\Phi^{T}V}/{N}$. Firstly, apart from the fixed $\theta_{0}$, we consider asymptotic properties of $N(\Phi^{T}\Phi)^{-1}$ and ${\Phi^{T}V}/{N}$ in Theorem \ref{thm:as and d convergence for terms of ls estimator}.

 In the following part, we adopt the concepts of almost sure convergence and convergence in distribution. We define that the random sequence $\{\xi_{N}\}\in\R^{d}$ converges almost surely to a random variable $\xi\in\R^{d}$ if $\Pr(\lim_{N\to\infty}\|\xi_{N}-\xi\|_{2}=0)=1$, which can be written as $\xi_{N}\overset{a.s.}\to \xi$ as $N\to\infty$. Define that the random sequence $\{\xi_{N}\}$ converges in distribution to a random variable $\xi$, if $\lim_{N\to\infty}\text{Pr}(\xi_{N}\leq x)=\text{Pr}(\xi\leq x)$ for every $x$ at which the limit distribution function $\text{Pr}(\xi\leq x)$ is continuous, where the map $x\mapsto \text{Pr}(\xi\leq x)$ denotes the distribution function of $\xi$. It can be written as $\xi_{N} \overset{d}\to \xi$.

\begin{theorem}\label{thm:as and d convergence for terms of ls estimator}
 Under Assumption \ref{asp:input}, \ref{asp:noise} and \ref{asp:independence between input and noise}, as $N\to\infty$, we have
			\begin{gather}
			\label{eq:almost sure convergence of PP_N}
			\frac{\Phi^{T}\Phi}{N}\overset{a.s.}\to \Sigma,\\
            \label{eq:almost sure convergence of PV_N}
            \frac{\Phi^{T}V}{N}\overset{a.s.}\to 0,\\
			\label{eq:almost sure convergence of VV_N}
			\frac{V^{T}V}{N}\overset{a.s.}\to \sigma^2,\\
			\left(\begin{array}{ccc}\sqrt{N}\left(\frac{\Phi^{T}\Phi}{N}-\Sigma\right), & \sqrt{N}(\frac{\Phi^{T}V}{N}),
            & \sqrt{N}(\frac{V^{T}V}{N}-\sigma^2) \end{array}\right)\nonumber\\
            \label{eq:joint convergence in distribution of PP, PV and VV}
			\overset{d}\to \left(\begin{array}{ccc}\Gamma, & \upsilon, &\rho \end{array}\right),
			\end{gather}
			where $\Gamma\in\R^{n\times n}$, $\upsilon\in\R^{n\times 1}$ and $\rho\in\R$ are jointly Gaussian distributed with zero mean and satisfy
			\begin{align}\label{eq:covariance of Gamma}
			C_{\Gamma}\triangleq& \E(\Gamma\otimes\Gamma)\nonumber\\
			=&\lim_{N\to\infty}N\E\left[\left(\frac{\Phi^{T}\Phi}{N}-\Sigma\right)\otimes \left(\frac{\Phi^{T}\Phi}{N}-\Sigma\right)\right],\\
			\label{eq:covariance of upsilon}
			\E(\upsilon\upsilon^{T})=&\sigma^2\Sigma,\\
			\label{eq:variance of rho}
			\E(\rho^2)=&\E[v^{4}(t)]-\sigma^4,\\
			\label{eq:correlatedness of upsilon and Gamma}
			\E(\upsilon\otimes \Gamma)=&0\in\R^{n^2\times n},\\
			\label{eq:correlatedness of upsilon and rho}
			\E(\rho\upsilon)=&0\in\R^{n\times 1},\\
			\label{eq:corelatedness of Gamma and rho}
			\E(\rho\Gamma)=&0\in\R^{n\times n}.
			\end{align}
Here $\otimes$ denotes the Kronecker product. 

Moreover, for $i,j=1,\cdots,n$, the $(i,j)$th element of $\Sigma$ can be represented as
\begin{align}\label{eq:def of Sigma element}
\left[\Sigma \right]_{i,j}=R_{u}(|i-j|);
\end{align}
for $i,j=1,\cdots,n^2$, the $(i,j)$th element of $C_{\Gamma}$ can be represented as
\begin{align}\label{eq:explicit representation of Expectation of Gamma_o_Gamma}
[C_{\Gamma}]_{i,j}
=&\left\{\E[e^{4}(t)]/\sigma_{e}^4-3 \right\} R_{u}(k)R_{u}(l)\nonumber\\
&+\sum_{\tau=-\infty}^{\infty}\left[R_{u}(\tau)R_{u}(\tau+k-l)+R_{u}(\tau+k)R_{u}(\tau-l)  \right],
\end{align}
where $R_{u}(\tau)$ is defined in \eqref{eq:definition of R_u_tau}, and
\begin{subequations}\label{eq:def of k and l}
	\begin{align}
	k=&\left|\lfloor{(i-1)/n}\rfloor-\lfloor{(j-1)/n}\rfloor\right|,\\
	l=&|i-j-\lfloor{(i-1)/n}\rfloor n+\lfloor{(j-1)/n}\rfloor n|.
	\end{align}
\end{subequations}
Here $\lfloor{\cdot}\rfloor$ denotes the floor operation, i.e. $\lfloor{x}\rfloor=\max\{\tilde{x}\in\Z|\tilde{x}\leq x\}$.
\end{theorem}

Then, combining \eqref{eq:LS estimator form3} with Theorem \ref{thm:as and d convergence for terms of ls estimator}, we can derive the almost sure convergence and the convergence in distribution of the LS estimator.
\begin{theorem}\label{thm:properties of ls estimate}
Under Assumption \ref{asp:input}, \ref{asp:noise} and \ref{asp:independence between input and noise}, as $N\to\infty$, we have
			\begin{gather}
			\label{eq:almost sure convergence of NPPinv}
			N(\Phi^{T}\Phi)^{-1}\overset{a.s.}\to \Sigma^{-1},\\
			\label{eq:almost sure convergence of ls estimate}
			\hat{\theta}^{\LS}\overset{a.s.}\to  \theta_{0}\\
            \label{eq:almost sure convergence of noise variance estimator}
            \widehat{\sigma^2}\overset{a.s.}\to \sigma^2\\
            \left(\sqrt{N}\left(N(\Phi^{T}\Phi)^{-1}-\Sigma^{-1}\right),  \sqrt{N}(\hat{\theta}^{\LS}-\theta_{0}),
             \sqrt{N}\left(\widehat{\sigma^2}-\sigma^2\right) \right)\nonumber\\
            \label{eq:joint convergence in distribution of ls estimate and other terms}
			\overset{d}\to \left(\begin{array}{ccc}-\Sigma^{-1}\Gamma\Sigma^{-1}, & \Sigma^{-1}\upsilon, &\rho \end{array}\right),
			\end{gather}
where $\Gamma$, $\upsilon$ and $\rho$ are defined as Theorem \ref{thm:as and d convergence for terms of ls estimator}, and
\begin{align}\label{eq:noise variance estimator}
\widehat{\sigma^2}=\frac{\|Y-\Phi\hat{\theta}^{\LS}\|_{2}^{2}}{N-n}.
\end{align}
\end{theorem}

Moreover, the boundedness in probability of building blocks and the LS estimator can also be derived according to Theorem \ref{thm:as and d convergence for terms of ls estimator} and \ref{thm:properties of ls estimate}, from which we can observe their convergence rates. For the nonzero sequence $\{a_{N}\}$, we let ${\xi}_{N}=O_{p}(a_{N})$ denote that $\{{\xi_{N}}/a_{N}\}$ is bounded in probability, which means that $\forall \epsilon>0$, $\exists L>0$ such that $\sup_{N}\text{Pr}(\|{\xi}_{N}/a_{N}\|_{2}>L)<\epsilon$.
\begin{theorem}\label{thm:boundedness in probability of ls estimate}
 Under Assumption \ref{asp:input}, \ref{asp:noise} and \ref{asp:independence between input and noise}, we have
            \begin{align}
			\label{eq:boundedness in probability of PP_inv}
			\Phi^{T}\Phi=&O_{p}(N),\\
            \label{eq:boundedness in probability of VV_N}
            V^{T}V=&O_{p}(N),\\
			\label{eq:boundedness in probability of PP_N}
			\frac{\Phi^{T}\Phi}{N}-\Sigma=& O_{p}(1/\sqrt{N}),\\
            \label{eq:boundedness in probability of PV}
			\frac{\Phi^{T}V}{N}=&O_{p}(1/\sqrt{N}),\\
            \label{eq:boundedness in probability of VV_N_sigma2}
            \frac{V^{T}V}{N}-\sigma^2=&O_{p}(1/\sqrt{N}),\\
			\label{eq:boundedness in probability of ls estimate}
			\hat{\theta}^{\LS}-\theta_{0}=&O_{p}(1/\sqrt{N}),\\
            \label{eq:boundedness in probability of noise variance estimator}
            \widehat{\sigma^2}-\sigma^2=&O_{p}(1/\sqrt{N}),
			\end{align}
where $\widehat{\sigma^2}$ is defined as \eqref{eq:noise variance estimator}.
\end{theorem}

\section{Boundedness of Moments of Several Terms}

In this section, we show the boundedness of moments of several terms.

\begin{assumption}\label{asp:bounded 8th moment of et and vt}
	The $8$th moments of $e(t)$ and $v(t)$ are both bounded.
\end{assumption}

\begin{assumption}\label{asp:bounded 16th moment of et and vt}
	The $16$th moments of $e(t)$ and $v(t)$ are both bounded.
\end{assumption}

\begin{theorem}\label{thm:boundedness of 4th moments of PVN}
	There exists $\widetilde{M}>0$, which is irrespective of $N$, such that under Assumptions \ref{asp:input}-\ref{asp:independence between input and noise}, we have
		\begin{align}
			\label{eq:boundedness of 4th moment of PVN}
			\E\left(\left\|\frac{\Phi^{T}V}{N}\right\|_{2}^4\right)\leq \frac{1}{N^2}\widetilde{M}.
		\end{align}
\end{theorem}

\begin{theorem}\label{thm:boundedness of 8th moment of PVN and 4th moment of VVN}
	There exists $\widetilde{M}>0$, which is irrespective of $N$, such that under Assumption \ref{asp:input}-\ref{asp:bounded 8th moment of et and vt}, we have
	\begin{align}\label{eq:boundedness of 8th moment of PVN}
		\E\left(\left\|\frac{\Phi^{T}V}{N}\right\|_{2}^8\right)\leq \frac{1}{N^4}\widetilde{M},\\
		\E\left(\left\|\frac{\Phi^{T}\Phi}{N}-\Sigma\right\|_{F}^{4} \right)\leq \frac{1}{N^2}\widetilde{M},\\
		\E\left(\left|\frac{V^{T}V}{N}-\sigma^2\right|^4 \right)\leq \frac{1}{N^2}\widetilde{M}.
	\end{align}
\end{theorem}

\begin{theorem}\label{thm:boundedness of 8th moment of PPN and VVN}
	There exists $\widetilde{M}>0$, which is irrespective of $N$, such that under Assumption \ref{asp:input}-\ref{asp:bounded 16th moment of et and vt}, we have
	\begin{align}\label{eq:boundedness of 8th moment of diff of PPN and Sigma}
		\E\left(\left\|\frac{\Phi^{T}\Phi}{N}-\Sigma\right\|_{F}^{8} \right)\leq \frac{1}{N^4}\widetilde{M},\\
		\label{eq:boundedness of 8th moment of diff of noise variance estimator}
		\E\left(\left|\frac{V^{T}V}{N}-\sigma^2\right|^8 \right)\leq \frac{1}{N^4}\widetilde{M}.
	\end{align}
\end{theorem}

\begin{remark}
	On the one hand,
	if Assumption \ref{asp:bounded 16th moment of et and vt} is true, then Assumption \ref{asp:bounded 8th moment of et and vt} must be true, because the boundedness of higher order moments always implies the boundedness of lower order moments. On the other hand
	if $e(t)$ and $v(t)$ are both Gaussian distributed, as long as their second moments are bounded, both Assumption \ref{asp:bounded 8th moment of et and vt} and \ref{asp:bounded 16th moment of et and vt} are true. 
\end{remark}

\section{Preliminary Results about Regularized Least Square Estimator for FIR Model}

To handle the ill-conditioned problem, one can introduce a regularization term in \eqref{eq:LS form1} to obtain the regularized least squares (RLS) estimator:
\begin{subequations}\label{eq:RLS estimator}
	\begin{align}
	\hat{\theta}^{\TR}=&\argmin_{\theta\in\R^{n}}\|Y-\Phi\theta\|_{2}^{2}+\sigma^2\theta^{T}P^{-1}\theta\\
	\label{eq:rls estimator form1}
	=&(\Phi^{T}\Phi+\sigma^2 P^{-1})^{-1}\Phi^{T}Y\\
	=&P\Phi^{T}Q^{-1}Y,
	\end{align}
\end{subequations}
where $P\in\R^{n\times n}$ is positive semidefinite, its $(i,j)$th element $P_{i,j}$ can be designed through a positive semidefinite \emph{kernel} $\kappa(i,j;\eta):\mathbb N\times \mathbb N\to \R$ with $\eta\in\Omega\subset\R^{p}$ being the hyper-parameter 
and thus $P$ is often called the kernel matrix, and
\begin{align}
Q=\Phi P\Phi^{T}+\sigma^2I_{N},
\end{align}
and $I_{N}$ denotes the $N$-dimensional identity matrix.

Moreover, we define that
\begin{align}
\label{eq:def of hat_S_inv}
\hat{S}(\eta)=&P(\eta)+\widehat{\sigma^2}(\Phi^{T}\Phi)^{-1}.
\end{align}

\begin{theorem}\label{thm:preliminary results about P and S}
	For the FIR model \eqref{eq:vector-matrix form linear regression model}, under Assumption \ref{asp:input}-\ref{asp:independence between input and noise}, we have the following results.
	\begin{enumerate}
		\item For any given $\eta\in\R^{p}$, we have
		\begin{align}
		\label{eq:almost sure convergence of sqrtN_NPPinv}
		\sqrt{N}(\Phi^{T}\Phi)^{-1}\overset{a.s.}\to& 0,\\
		\label{eq:almost sure convergence of S_inv}
		\hat{S}(\eta)^{-1}\overset{a.s.}\to & P(\eta)^{-1},\\
		\label{eq:almost sure convergence of sqrtN_Sinv_Pinv}
		\sqrt{N}(\hat{S}(\eta)^{-1}-P(\eta)^{-1})\overset{a.s.}\to& 0,\\
		\label{eq:almost sure convergence of 1st order derivatives of S_inv and P_inv}
		\frac{\partial \hat{S}(\eta)^{-1}}{\partial \eta_{k}} \overset{a.s.}\to& \frac{\partial P(\eta)^{-1}}{\partial \eta_{k}},\\
		\label{eq:almost sure convergence of sqrtN difference of 1st order derivatives of S_inv and P_inv}
		\sqrt{N}\left(\frac{\partial \hat{S}(\eta)^{-1}}{\partial \eta_{k}}-\frac{\partial P(\eta)^{-1}}{\partial \eta_{k}}\right) \overset{a.s.}\to& 0,
		\end{align}
		where $\eta_{k}$ denotes the $k$th element of $\eta$ and $k=1,\cdots,p$.
		
		\item Suppose that as $N\to\infty$, $\hat{\eta}_{N}\overset{a.s.}\to \eta^{*}$. If $P(\eta)$ is continuous for every $\eta\in\widetilde{\Omega}$ and there exists a compact set $\widetilde{\Omega}$ containing $\eta^{*}$ such that $\|\hat{S}^{-1}\|_{F}<\|P^{-1}\|_{F}$ is bounded, we have as $N\to\infty$,
		\begin{align}
		\label{eq:almost sure convergence of S_inv to P_inv at convergent estimate}
		\hat{S}(\hat{\eta}_{N})^{-1} \overset{a.s.}\to P(\eta^{*})^{-1}.
		\end{align}
		
		\item Suppose that as $N\to\infty$, $\hat{\eta}_{N}\overset{a.s.}\to \eta^{*}$. If $P(\eta)$ is differentiable for every $\eta\in\Omega$, then we have
		\begin{align}
		\label{eq:diff of S_inv and P_inv for general P}
		&\hat{S}(\hat{\eta}_{N})^{-1}-P(\eta^{*})^{-1}\nonumber\\
		=&-\hat{S}(\hat{\eta}_{N})^{-1}
		\left[\sum_{k=1}^{p}\left.\frac{\partial P(\eta)}{\partial \eta_{k}}\right|_{\eta=\tilde{\eta}_{N}}e_{k}^{T}(\hat{\eta}_{N}-\eta^{*})\right] P(\eta^{*})^{-1}\nonumber\\
		&-\widehat{\sigma^2}\hat{S}(\hat{\eta}_{N})^{-1}(\Phi^{T}\Phi)^{-1}P(\eta^{*})^{-1},
		\end{align}
		where $e_{k}\in\R^{p}$ denotes a column vector with $k$th element being one and others zero, and $\tilde{\eta}_{N}$ belongs to a neighborhood of $\eta^{*}$ with radius $\|\hat{\eta}_{N}-\eta^{*}\|_{2}$.
		In particular, if we consider $P=\eta I_{n}$ with $\eta>0$, we have
		\begin{align}
		\label{eq:equivalent form of diff of S_inv and P_inv}
		&\hat{S}(\hat{\eta}_{N})^{-1}-P(\eta^{*})^{-1}\nonumber\\
		=&-(\hat{\eta}_{N}-\eta^{*})\hat{S}(\hat{\eta}_{N})^{-1}P(\eta^{*})^{-1}\nonumber\\
		-&\widehat{\sigma^2}\hat{S}(\hat{\eta}_{N})^{-1}(\Phi^{T}\Phi)^{-1}P(\eta^{*})^{-1}.
		\end{align}
		
		\item We have
		\begin{align}\label{eq:difference of S_inv and P_inv}
		\hat{S}(\eta)^{-1}-P(\eta)^{-1}=-\frac{1}{N}\widehat{\sigma^2} \hat{S}(\eta)^{-1}N(\Phi^{T}\Phi)^{-1}P(\eta)^{-1}.
		\end{align}
		For $k,l=1,\cdots,p$, we have
		\begin{subequations}\label{eq:1st derivatives of matrix inverse}
			\begin{align}\label{eq:equivalent form of 1st derivetive of P inverse}
			\frac{\partial P(\eta)^{-1}}{\partial \eta_{k}} =& -P^{-1}\frac{\partial P(\eta)}{\partial \eta_{k}}P^{-1},\\
			\frac{\partial \hat{S}(\eta)^{-1}}{\partial \eta_{k}} =& -\hat{S}(\eta)^{-1}\frac{\partial P(\eta)}{\partial \eta_{k}}\hat{S}(\eta)^{-1},
			\end{align}
		\end{subequations}
		\begin{subequations}\label{eq:2nd derivatives of matrix inverse}
			\begin{align}
			\frac{\partial^2 P(\eta)^{-1}}{\partial\eta_{k}\partial\eta_{l}}
			=&P(\eta)^{-1}\frac{\partial P(\eta)}{\partial \eta_{l}}P(\eta)^{-1}\frac{\partial P(\eta)}{\partial \eta_{k}}P(\eta)^{-1}\nonumber\\
			&-P(\eta)^{-1}\frac{\partial^2 P(\eta)}{\partial\eta_{k}\partial\eta_{l}}P(\eta)^{-1}\nonumber\\
			&+P(\eta)^{-1}\frac{\partial P(\eta)}{\partial \eta_{k}}P(\eta)^{-1}\frac{\partial P(\eta)}{\partial\eta_{l}}P(\eta)^{-1},\\
			\frac{\partial^2 \hat{S}(\eta)^{-1}}{\partial\eta_{k}\partial\eta_{l}}
			=&\hat{S}(\eta)^{-1}\frac{\partial P(\eta)}{\partial \eta_{l}}\hat{S}(\eta)^{-1}\frac{\partial P(\eta)}{\partial \eta_{k}}\hat{S}(\eta)^{-1}\nonumber\\
			&-\hat{S}(\eta)^{-1}\frac{\partial^2 P(\eta)}{\partial\eta_{k}\partial\eta_{l}}\hat{S}(\eta)^{-1}\nonumber\\
			&+\hat{S}(\eta)^{-1}\frac{\partial P(\eta)}{\partial \eta_{k}}\hat{S}(\eta)^{-1}\frac{\partial P(\eta)}{\partial\eta_{l}}\hat{S}(\eta)^{-1}.
			\end{align}
		\end{subequations}
	\end{enumerate}
\end{theorem}

\section{Conclusion}

For the FIR model with filtered white noise inputs, we mainly consider asymptotic properties of the LS estimator in terms of: the almost sure convergence, the convergence in distribution and the boundedness in probability. 
Moreover, some preliminary results of the RLS estimator are also contained. These properties help us have a better understanding of the LS and RLS estimators. In addition, we can use these tools to investigate asymptotic properties of other estimators, such as hyper-parameter estimators.

\def\thesectiondis{\thesection.}                   
\def\thesubsectiondis{\thesection.\arabic{subsection}.}          
\def\thesubsubsectiondis{\thesubsection.\arabic{subsubsection}.}

\setcounter{subsection}{0}

\renewcommand{\thesection}{A}
\setcounter{theorem}{0}
\renewcommand{\thelemma}{A.\arabic{lemma}}

\renewcommand{\theequation}{A.\arabic{equation}}
\setcounter{equation}{0}

\renewcommand{\thesubsection}{\thesection.\arabic{subsection}}

\section*{Appendix A}\label{sec:Appendix A}

All proofs of theorems and corollaries are included in Appendix A, and all required lemmas and their corresponding proofs are shown in Appendix B.

\subsection{Proof of Theorem \ref{thm:as and d convergence for terms of ls estimator}}

\begin{enumerate}
    \item \underline{Proof of \eqref{eq:almost sure convergence of PP_N}}

    It can be seen that for $i,j=1,\cdots,n$, the $(i,j)$th element of $\Phi^{T}\Phi/N$ is
    \begin{align}
    \left[\frac{\Phi^{T}\Phi}{N}\right]_{i,j}=\sum_{t=1-i}^{N-i}u(t)u(t+i-j).
    \end{align}
	Using Lemma \ref{lemma:ergodic theory} with
	\begin{align}
	x(t)=u(t),\ m(t)=0,\ R_{s}(\tau)=R_{u}(\tau),
	\end{align}
	we obtain that as $N\to\infty$,
	\begin{align}\label{eq:term uu}
	\frac{1}{N}\sum_{t=1}^{N}u(t)u(t-\tau) \overset{a.s.}\to R_{u}(\tau).
	\end{align}
	Applying \eqref{eq:term uu} to each element of $\Phi^{T}\Phi/N$, since the almost sure convergence of a matrix is equivalent to that of its all elements, it completes the proof of \eqref{eq:almost sure convergence of PP_N}.

    \item \underline{Proof of \eqref{eq:almost sure convergence of PV_N}}

    Since
	\begin{align}
	\frac{\Phi^{T}V}{N}=\left[\begin{array}{c}\frac{1}{N}\sum_{t=0}^{N-1}u(t)v(t+1)\\ \frac{1}{N}\sum_{t=-1}^{N-2}u(t)v(t+2)\\ \vdots \\ \frac{1}{N}\sum_{t=1-n}^{N-n}u(t)v(t+n)\end{array}\right],
	\end{align}
	we can consider the almost sure convergence of the general form of elements. It means that as long as $\frac{1}{N}\sum_{t=1}^{N}u(t)v(t+i)\overset{a.s.}\to 0$ with $i=1,\cdots,n$ when $N\to\infty$, ${\Phi^{T}V}/{N}\overset{a.s.}\to 0$ holds.
	
	Recall that
	\begin{align}
	\frac{1}{N}\sum_{t=1}^{N}u(t)v(t+i)=\frac{1}{N}\sum_{t=1}^{N}\sum_{k=0}^{\infty}h(k)e(t-k)v(t+i).
	\end{align}
	
	Firstly, using Lemma \ref{lem:Boundedness of Expectation} with
	\begin{align}
	w_{1}(t)=u(t),\ w_{3}(t)=v(t),\ r=1,
	\end{align}
	we obtain that
	\begin{align}\label{eq:middle result for bounded expectation of PV}
	\E\left[\sum_{t=1}^{N}u(t)v(t+i)\right]^2\leq C_{1}N,
	\end{align}
	where $C_{1}=\left[\sum_{k=0}^{\infty}|h(k)|\right]^2\sigma^2\sigma_{e}^2$.
	
	Then we apply Lemma \ref{lemma:Almost Sure Convergence for Statistics with Bounded Expectation} with
	\begin{align}
	\tilde{s}(t)=u(t)v(t+i),\ C_{s}=C_{1},
	\end{align}
	hence as $N\to\infty$, it yields that
	\begin{align}
	\frac{1}{N}\sum_{t=1}^{N}u(t)v(t+i)\overset{a.s.}\to 0,
	\end{align}
	which can be extended to each element of $(\Phi^{T}V)/N$. It means that as $N\to\infty$, we have \eqref{eq:almost sure convergence of PV_N}.

	\item \underline{Proof of \eqref{eq:almost sure convergence of VV_N}}
	
	For
	\begin{align}
	\frac{V^{T}V}{N}=\frac{1}{N}\sum_{t=1}^{N}v^2(t),
	\end{align}
	since for $t=1,2,\cdots,N$,
	\begin{align}
	\E[v^2(t)]=\E|v^2(t)|=\sigma^2<\infty,
	\end{align}
	we can derive \eqref{eq:almost sure convergence of VV_N} from Lemma \ref{lem:kolmogorov strong law}.
	
	\item \underline{Proof of \eqref{eq:joint convergence in distribution of PP, PV and VV}, \eqref{eq:covariance of Gamma}, \eqref{eq:covariance of upsilon}, \eqref{eq:variance of rho}, \eqref{eq:correlatedness of upsilon and Gamma}, \eqref{eq:correlatedness of upsilon and rho} and \eqref{eq:corelatedness of Gamma and rho}}
	
	It is equivalent to derive the convergence in distribution of the following vector
	\begin{align}
	S_{N}=\left[\begin{array}{c} \text{vec}(\sqrt{N}(\Phi^{T}\Phi/N-\Sigma))\\ \sqrt{N}\Phi^{T}V/N  \\ \sqrt{N}(V^{T}V/N-\sigma^2) \end{array}\right].
	\end{align}
	The $k$th element of $S_{N}\in\R^{n^2+n+1}$ is
	\begin{align}
	\left\{\begin{array}{ll}\frac{1}{\sqrt{N}}\sum_{t=1}^{N}\left[u(t-i)u(t-j)\right. & \text{if}\ 1\leq k\leq n^2, \\
     \left.-R_{u}(|i-j|)\right], & \\
	\frac{1}{\sqrt{N}}\sum_{t=1}^{N}u(t-k+n^2)v(t), & \text{if}\ n^2+1\leq k\leq n^2+n,\\
	\frac{1}{\sqrt{N}}\sum_{t=1}^{N}v^2(t), & \text{if}\ k=n^2+n+1, \end{array}\right.
	\end{align}
	where $k=1,2,\cdots,n^2+n+1$ and
	\begin{subequations}\label{eq:definition of i and j}
		\begin{align}
		i=&\lfloor{(k-1)/n}\rfloor+1\\
		j=&k-\lfloor{(k-1)/n}\rfloor n.
		\end{align}
	\end{subequations}
	Here $\text{vec}(\cdot)$ denotes the vectorization of a matrix and $\lfloor{\cdot}\rfloor$ denotes the floor operation, i.e. $\lfloor{x}\rfloor=\max\{\tilde{x}\in\Z|\tilde{x}\leq x\}$.
	
	For a positive integer $M$, let $u(t)=u^{M}(t)+\tilde{u}^{M}(t)$ with
	\begin{align}
	u^{M}(t)=&\sum_{k=0}^{M}h(k)e(t-k)\\
	\tilde{u}^{M}(t)=&\sum_{k=M+1}^{\infty}h(k)e(t-k).
	\end{align}
	Thus, we have
	\begin{align}\label{eq:split of S_N}
	S_{N}=Z_{M}(N)+X_{M}(N),
	\end{align}
	where
	\begin{align}
	Z_{M}(N)=&\sum_{t=1}^{N}\kappa_{M,N}(t)\\
	X_{M}(N)=&\frac{1}{\sqrt{N}}\tilde{\kappa}_{M}(N)
	\end{align}
	with their structure satisfying
	\begin{align}
	&[\kappa_{M,N}(t)]_{k}\nonumber\\
   =&\left\{\begin{array}{ll}\frac{1}{\sqrt{N}}\left\{u^{M}(t-i)u^{M}(t-j)\right. & \text{if}\ 1\leq k\leq n^2,\\ \left.-\E\left[u^{M}(t-i)u^{M}(t-j)\right]\right\}, & \\
	\frac{1}{\sqrt{N}}u^{M}(t-k+n^2)v(t), & \text{if}\ n^2+1\leq k\leq n^2+n,\\
	\frac{1}{\sqrt{N}}\left[v^2(t)-\sigma^2\right], & \text{if}\ k=n^2+n+1, \end{array}\right.\\
	&[\tilde{\kappa}_{M}(N)]_{k}\nonumber\\
	=&\left\{\begin{array}{ll}\sum_{t=1}^{N}\left\{\tilde{u}^{M}(t-i)u(t-j)\right. & \text{if}\ 1\leq k\leq n^2,\\
           +u^{M}(t-i)\tilde{u}^{M}(t-j) & \\
           -\E\left[\tilde{u}^{M}(t-i)u(t-j)\right. &\\
           +\left.\left.u^{M}(t-i)\tilde{u}^{M}(t-j)\right] \right\}, &\\
	\sum_{t=1}^{N}\tilde{u}^{M}(t-k+n^2)v(t), & \text{if}\ n^2+1\leq k\leq n^2+n,\\
	0, & \text{if}\ k=n^2+n+1, \end{array}\right.
	\end{align}
	where $i$ and $j$ are defined as \eqref{eq:definition of i and j}. In the following part, we shall apply Lemma \ref{lem:Limiting theorem of sum} for $S_{N}$.
	
	\begin{itemize}
		\item [$-$] Firstly, we will show that as $N\to\infty$,
		\begin{align}
		Z_{M}(N)\overset{d}\to& \mathcal{N}(0,Q_{M}),\\
		Q_{M}=&\lim_{N\to\infty}\E[Z_{M}(N)Z_{M}(N)^{T}].
		\end{align}
		
		It can be seen that $\E(\kappa_{M,N}(t))=0$.
		And $\kappa_{M,N}(t)$ is $(M+n-1)$-dependent, which means that
		\begin{align}
		\{\kappa_{M,N}(1),\cdots,\kappa_{M,N}(s)\}\ \text{and}\ \{\kappa_{M,N}(t),\kappa_{M,N}(t+1),\cdots\}
		\end{align}
		are independent if $t-s>M+n-1$, which can be derived from $\max_{i,j}|i-j|=n-1$.
		
		Then for
		\begin{align}
		\E|u^{M}(t)|^{4+\delta}=&\E\left|\sum_{k=0}^{M}h(k)e(t-k)\right|^{4+\delta}\nonumber\\
		\leq& 2^{(3+\delta)M}\sum_{k=0}^{M}\E|h(k)e(t-k)|^{4+\delta}\nonumber\\
		&[\text{using}\ \text{Lemma}\ \ref{lem:C_r inequality}]\nonumber\\
		\leq& 2^{(3+\delta)M}\sup_{t}\E|e(t)|^{4+\delta}\left[\sum_{k=0}^{M}|h(k)|\right]^{4+\delta}\nonumber\\
		\leq& 2^{(3+\delta)M}\sup_{t}\E|e(t)|^{4+\delta}\left[\sum_{k=0}^{\infty}|h(k)|\right]^{4+\delta},
		\end{align}
		since $\sum_{k=0}^{\infty}|h(k)|<\infty$ and $\E|e(t)|^{4+\delta}<\infty$ for some $\delta>0$, we know that $u^{M}(t)$ also has bounded moments of order $4+\delta$ for some $\delta>0$. Then we apply Lemma \ref{lemma:bounded result} with
		\begin{itemize}
			\item [$\star$] when $1\leq k\leq n^2$, $\beta_{k}=u^{M}(t-i)$ and $\gamma_{k}=u^{M}(t-j)$;
			\item [$\star$] when $n^2+1\leq k \leq n^2+n$, $\beta_{k}=u^{M}(t-k+n^2)$ and $\gamma_{k}=v(t)$;
			\item [$\star$] when $k=n^2+n+1$, $\beta_{k}=v(t)$ and $\gamma_{k}=v(t)$
		\end{itemize}
	    to obtain that there exists a constant $C_{1}>0$ such that
	    \begin{align}
	    \E\|\kappa_{M,N}(t)\|_{2}^{2+\delta}\leq \frac{1}{N}N^{-\delta/2}C_{1}.
	    \end{align}
	    It follows that
	    \begin{align}
	    &{\lim\sup}_{N\to\infty}\sum_{t=1}^{N}\E\|\kappa_{M,N}(t)\|_{2}^2<\infty\\
	    &\lim_{N\to\infty}\sum_{t=1}^{N}\E\|\kappa_{M,N}(t)\|_{2}^{2+\delta}=0\ \text{for}\ \text{some}\ \delta>0.
	    \end{align}
	
	    Thus, we can apply Lemma \ref{lem:CLT for M-dependent sequence} with $x_{N}(t)=\kappa_{M,N}(t)$ to yield as $N\to\infty$,
	    \begin{align}
	    Z_{M}(N)\overset{d}\to& \mathcal{N}(0,Q_{M}),\\
	    Q_{M}=&\lim_{N\to\infty}\E[Z_{M}(N)Z_{M}(N)^{T}].
	    \end{align}
	
	    \item [$-$] The second step is to show that there exists a constant ${C}_{M}>0$ such that
		\begin{align}
		\E\|X_{M}(N)\|_{2}^2\leq &{C}_{M},\\
		\lim_{M\to\infty}{C}_{M}=&0.
		\end{align}
		
        Note that when $k=n^2+n+1$, the $k$th element of $\tilde{\kappa}_{M}(N)$ is zero. Then for
		\begin{align}
		&\E\|\tilde{\kappa}_{M}(N)\|_{2}^{2}\nonumber\\
		=&\sum_{k=1}^{n^2}\E\left\{\sum_{t=1}^{N}\left[\tilde{u}^{M}(t-i)u(t-j)+u^{M}(t-i)\tilde{u}^{M}(t-j)\right.\right.\nonumber\\
		&-\left.\left.\E[\tilde{u}^{M}(t-i)u(t-j)+u^{M}(t-i)\tilde{u}^{M}(t-j)]\right] \right\}^2\nonumber\\
		&+\sum_{k=n^2+1}^{n^2+n}\E\left\{\sum_{t=1}^{N}\tilde{u}^{M}(t-k+n^2)v(t) \right\}^2,
		\end{align}
		if we apply the adjusted version of Lemma \ref{lem:Boundedness of Expectation}, there must exist a constant $C_{3}$ such that
		\begin{align}
		\E\|\tilde{\kappa}_{M}(N)\|_{2}^{2}\leq C_{3}\left[\sum_{k=M+1}^{\infty}|h(k)|\right]^{2}N,
		\end{align}
		which gives that
		\begin{align}
		\E\|X_{M}(N)\|_{2}^2\leq C_{M}=C_{3}\left[\sum_{k=M+1}^{\infty}|h(k)|\right]^{2}.
		\end{align}
		In addition, since $\sum_{k=0}^{\infty}|h(k)|$ is always finite, it can be seen that $\sum_{k=0}^{N}|h(k)|$ is convergent, leading to
		\begin{align}
		\lim_{M\to\infty}C_{M}=\lim_{M\to\infty}C_{3}\left[\sum_{k=M+1}^{\infty}|h(k)|\right]^{2}=0.
		\end{align}	
	
	    \item [$-$] Combining two steps above with \eqref{eq:split of S_N}, we can use Lemma \ref{lem:Limiting theorem of sum} to know that as $N\to\infty$,
	    \begin{align}
        \label{eq:convergence in distribution of S_N}
	    S_{N}\overset{d}\to& (0,Q_{s})\\
        \label{eq:def of coveriance matrix}
	    Q_{s}=&\lim_{M\to\infty}Q_{M}=\lim_{N\to\infty}\E[S_{N}S_{N}^{T}],
	    \end{align}
	    where the last step comes from $\lim_{M\to\infty}\lim_{N\to\infty}\E\|X_{M}(N)\|_{2}^2=0$.
	\end{itemize}

    Hence, from \eqref{eq:convergence in distribution of S_N}, we can equivalently derive \eqref{eq:joint convergence in distribution of PP, PV and VV} as $N\to\infty$. Moreover, from \eqref{eq:def of coveriance matrix}, we have \eqref{eq:covariance of Gamma} and the following results.
    \begin{itemize}
        \item [$\circ$] For \eqref{eq:covariance of upsilon}, we have
        \begin{align}
        \E(\upsilon\upsilon^{T})=&\lim_{N\to\infty}\frac{1}{N}\E(\Phi^{T}VV^{T}\Phi)\nonumber\\
	    =&\lim_{N\to\infty}\sigma^2\frac{\Phi^{T}\Phi}{N}\nonumber\\
	    =&\sigma^2\Sigma.
        \end{align}

        \item [$\circ$] For \eqref{eq:correlatedness of upsilon and Gamma}, we have
        \begin{align}\label{eq:expectation of upsilon and Gamma}
        \E(\upsilon\otimes \Gamma)=\lim_{N\to\infty}N\E\left\{[\Phi^{T}V/N]\otimes [(\Phi^{T}\Phi)/N-\Sigma]\right\}.
        \end{align}
    Then for $k_{1},k_{2},l_{2}=1,\cdots,n$, each element of \eqref{eq:expectation of upsilon and Gamma} can be represented as follows,
	\begin{align}
	&\E\left\{ \left[\frac{\Phi^{T}V}{N}\right]_{k_{1}}\left[ \frac{\Phi^{T}\Phi}{N}-\Sigma\right]_{k_{2},l_{2}} \right\}\nonumber\\
	=&\E\left\{ \frac{1}{N^2}\sum_{t_{1}=1}^{N}u(t_{1}-k_{1})v(t_{1})\right.\nonumber\\
	&\left.\sum_{t_{2}=1}^{N}\left[ u(t_{2}-k_{2})u(t_{2}-l_{2})-R_{u}(|k_{2}-l_{2}|)\right]  \right\}\nonumber\\
	=&\frac{1}{N^2}\sum_{t_{1}=1}^{N}\sum_{t_{2}=1}^{N} \E[v(t_{1})]\nonumber\\
	&\E\left[u(t_{1}-k_{1})\left(u(t_{2}-k_{2})u(t_{2}-l_{2})-R_{u}(|k_{2}-l_{2}|)\right)\right],
	\end{align}
	where the last step comes from the mutual independence of $\{v(t)\}_{t=1}^{N}$ and $\{u(t)\}_{t=1-n}^{N-1}$. Since
	\begin{align}
	&\E(v(t_{1}))=0,\\
	&\E\left[u(t_{1}-k_{1})\left(u(t_{2}-k_{2})u(t_{2}-l_{2})-R_{u}(|k_{2}-l_{2}|)\right)\right]\nonumber\\
	=&\E\left[u(t_{1}-k_{1})u(t_{2}-k_{2})u(t_{2}-l_{2})\right]\nonumber\\
	&-\E[u(t_{1}-k_{1})]R_{u}(|k_{2}-l_{2}|)\nonumber\\
	=&\E\left[\sum_{m_{1}=0}^{\infty}\sum_{m_{2}=0}^{\infty}\sum_{m_{3}=0}^{\infty} h(m_{1})h(m_{2})h(m_{3})\right.\nonumber\\
	&\left.e(t_{1}-k_{1}-m_{1})e(t_{2}-k_{2}-m_{2})e(t_{2}-l_{2}-m_{3})\right]\nonumber\\
	\leq& \sup_{t}\E\left[ e^3(t)\right]\left[\sum_{m=0}^{\infty}|h(m)|\right]^3<\infty,
	\end{align}
    we obtain \eqref{eq:correlatedness of upsilon and Gamma}.

    	\item [$\circ$] For \eqref{eq:variance of rho}, we have
    	\begin{align}
    	&\E(\rho^2)\nonumber\\
    	=&\lim_{N\to\infty}\frac{1}{N}\E\left[\sum_{t=1}^{N}\left(v^2(t)-\sigma^2\right)\right]^2\nonumber\\
    	=&\lim_{N\to\infty}\frac{1}{N}\sum_{t=1}^{N}\sum_{j=1}^{N}\E\left[ (v^2(t)-\sigma^2)(v^2(j)-\sigma^2)\right]\nonumber\\
    	=&\lim_{N\to\infty}\frac{1}{N}\sum_{t=1}^{N}\left\{ \E\left[v^4(t)\right]-2\sigma^2\E\left[v^2(t)\right] +\sigma^4 \right\}\nonumber\\
    	&[\text{since}\ v(t)\ \text{and}\ v(j)\ \text{are}\ \text{independent}\ \text{when}\ t\not=j]\nonumber\\
    	=&\E\left[v^4(t)\right] - \sigma^4.
    	\end{align}
    	
    	\item [$\circ$] For \eqref{eq:correlatedness of upsilon and rho}, we have
    	\begin{align}
    	\E(\rho\upsilon)=&\lim_{N\to\infty}\E\left[\Phi^{T}V\left(V^{T}V/N -\sigma^2 \right)\right].
    	\end{align}
    	The $k$th element of $\E(\rho\upsilon)$ with $k=1,2,\cdots,n$ can be represented as
    	\begin{align}
    	&\left[\E(\rho\upsilon) \right]_{k}\nonumber\\
    	=&\lim_{N\to\infty}\E\left\{ \left[\sum_{t=1}^{N}u(t-k)v(t) \right] \left[\frac{1}{N}\sum_{j=1}^{N}v^2(j) -\sigma^2 \right] \right\}\nonumber\\
    	=&\lim_{N\to\infty}\E\left\{ \sum_{t=1}^{N}\sum_{j=1}^{N}u(t-k)v(t)v^2(j) \right\} \nonumber\\
    	=&\lim_{N\to\infty}\sum_{t=1}^{N}\E\left[ u(t-k)v^3(t)\right]\nonumber\\
    	=&0.
    	\end{align}
    	
    	\item [$\circ$] For \eqref{eq:corelatedness of Gamma and rho}, we have
    	\begin{align}
    	\E(\rho\Gamma)
    	=&\lim_{N\to\infty}N\E\left[ \left(\frac{V^{T}V}{N} - \sigma^2\right) \left( \frac{\Phi^{T}\Phi}{N}-\Sigma \right) \right]\nonumber\\
    	=&\lim_{N\to\infty}N\E \left(\frac{V^{T}V}{N} - \sigma^2\right) \E\left( \frac{\Phi^{T}\Phi}{N}-\Sigma \right)\nonumber\\
    	=&0,
    	\end{align}
    	where the last second step derives from the independence between $\Phi$ and $V$.
    \end{itemize}
	\item \underline{Proof of \eqref{eq:def of Sigma element}}
	
	Combining \eqref{eq:zero mean of ut}, \eqref{eq:definition of R_u_tau} and \eqref{eq:def of Sigma}, for $i,j=1,\cdots,n$, it can be derived that 
	\begin{align}
	\left[ \Sigma \right]_{i,j}=\E\left[u(i-1)u(j-1)\right]=R_{u}(|i-j|).
	\end{align}
	
	\item \underline{Proof of \eqref{eq:explicit representation of Expectation of Gamma_o_Gamma}}
	
	For $i,j=1,\cdots,n^2$, the $(i,j)$th element of $C_{\Gamma}$  is
	\begin{align}
	&[C_{\Gamma}]_{i,j}\nonumber\\
	=&\lim_{N\to\infty}N\E\left[ \left( \frac{\Phi^{T}\Phi}{N}-\Sigma\right)_{k_{1},l_{1}} \left(  \frac{\Phi^{T}\Phi}{N}-\Sigma\right)_{k_{2},l_{2}}\right]\nonumber\\
	=&\lim_{N\to\infty}N\E\left\{ \left[ \frac{1}{N}\sum_{t_{1}=1}^{N}u(t_{1}-k_{1})u(t_{1}-l_{1})-R_{u}(|k_{1}-l_{1}|)\right]\right.\nonumber\\ &\left[\frac{1}{N}\sum_{t_{2}=1}^{N}u(t_{2}-k_{2})u(t_{2}-l_{2})-R_{u}(|k_{2}-l_{2}|)\right]\Bigg\}\nonumber\\
	=&\lim_{N\to\infty}N\E\left\{ \left[\frac{1}{N}\sum_{t_{1}=1}^{N}u(t_{1})u(t_{1}+|k_{1}-l_{1}|)-R_{u}(|k_{1}-l_{1}|)\right]\right.\nonumber\\ &\left[\frac{1}{N}\sum_{t_{2}=1}^{N}u(t_{2})u(t_{2}+|k_{2}-l_{2}|)-R_{u}(|k_{2}-l_{2}|)\right]\Bigg\},
	\end{align}
	where $k_{1}$, $k_{2}$, $l_{1}$ and $l_{2}$ satisfy 
	\begin{align}\label{eq:relationships between subscripts}
	k_{1}=&\lfloor{(i-1)/n}\rfloor+1\\
	l_{1}=&\lfloor{(j-1)/n}\rfloor+1\\
	k_{2}=&i-\lfloor{(i-1)/n}\rfloor n\\
	l_{2}=&j-\lfloor{(j-1)/n}\rfloor n.
	\end{align}
	
	For convenience, define
	\begin{subequations}\label{eq:def of p,q}
		\begin{align}
		k=&|k_{1}-l_{1}|,\\
		l=&|k_{2}-l_{2}|.
		\end{align}
	\end{subequations}
	Then applying Lemma \ref{lemma:properties of covariance function of filtered white noise},
	we can derive \eqref{eq:explicit representation of Expectation of Gamma_o_Gamma}.
	
\end{enumerate}

\subsection{Proof of Theorem \ref{thm:properties of ls estimate}}

\begin{enumerate}
	
\item \underline{Proof of \eqref{eq:almost sure convergence of NPPinv}}

Applying Lemma \ref{lemma:continuous mapping theorem for vector case} to \eqref{eq:almost sure convergence of PP_N}, we can obtain \eqref{eq:almost sure convergence of NPPinv}.

\item \underline{Proof of \eqref{eq:almost sure convergence of ls estimate}}

Combining \eqref{eq:almost sure convergence of PP_N} with \eqref{eq:almost sure convergence of PV_N}, we can apply Lemma \ref{lemma:Relationships between modes of convergence} and \ref{lemma:Slutsky Theorem} to \eqref{eq:LS estimator form3} to obtain \eqref{eq:almost sure convergence of ls estimate}.

\item \underline{Proof of \eqref{eq:almost sure convergence of noise variance estimator}}

We have
\begin{align}
&\frac{\|Y-\Phi\hat{\theta}^{\LS}\|_{2}^{2}}{N-n}\nonumber\\
=&\frac{\|\Phi(\theta_{0}-\hat{\theta}^{\LS})+V\|_{2}^{2}}{N-n}\nonumber\\
\label{eq:middle result of noise variance estimator}
=&-\frac{N}{N-n}(\hat{\theta}^{\LS}-\theta_{0})^{T}\frac{\Phi^{T}V}{N}+\frac{N}{N-n}\frac{V^{T}V}{N}.
\end{align}
Applying \eqref{eq:almost sure convergence of PP_N}, \eqref{eq:almost sure convergence of PV_N}, \eqref{eq:almost sure convergence of VV_N}, \eqref{eq:almost sure convergence of ls estimate} and Lemma \ref{lemma:continuous mapping theorem for vector case}, it leads to \eqref{eq:almost sure convergence of noise variance estimator}.

\item \underline{Proof of \eqref{eq:joint convergence in distribution of ls estimate and other terms}}

Firstly, based on  \eqref{eq:LS estimator form3}, \eqref{eq:almost sure convergence of PP_N} and Lemma \ref{lemma:Relationships between modes of convergence} and \ref{lemma:continuous mapping theorem for vector case}, we can obtain as $N\to\infty$,
\begin{align}\label{eq:convergence in distribution of ls estimate}
\sqrt{N}(\hat{\theta}^{\LS}-\theta_{0})\overset{d}\to \Sigma^{-1}\upsilon.
\end{align}

Then, we consider
\begin{align}
&\sqrt{N}\left( \widehat{\sigma^2}-\sigma^2\right)\nonumber\\
=&\sqrt{N}\left( \frac{\|Y-\Phi\hat{\theta}^{\LS}\|_{2}^{2}}{N-n}-\sigma^2 \right)\nonumber\\
=&-\frac{N}{N-n}\sqrt{N}(\hat{\theta}^{\LS}-\theta_{0})^{T}\frac{\Phi^{T}V}{N}\nonumber\\
&+\frac{N}{N-n}\sqrt{N}\left(\frac{V^{T}V}{N}-\sigma^2\right)+\frac{n\sigma^2}{N-n}.
\end{align}
Due to \eqref{eq:convergence in distribution of ls estimate}, \eqref{eq:almost sure convergence of PP_N}, \eqref{eq:almost sure convergence of PV_N} and \eqref{eq:joint convergence in distribution of PP, PV and VV} and Lemma \ref{lemma:Relationships between modes of convergence} and \ref{lemma:continuous mapping theorem for vector case}, it follows that as $N\to\infty$,
\begin{align}\label{eq:middle result for convergence in distribution of noise variance estimator}
&-\frac{N}{N-n}\sqrt{N}(\hat{\theta}^{\LS}-\theta_{0})^{T}\frac{\Phi^{T}V}{N}+\sqrt{N}n\sigma^2/(N-n)
\overset{d}\to 0.
\end{align}

Finally, combining \eqref{eq:convergence in distribution of ls estimate} and \eqref{eq:middle result of noise variance estimator}, we can apply \eqref{eq:middle result for convergence in distribution of noise variance estimator} and Lemma \ref{lemma:continuous mapping theorem for vector case} to derive \eqref{eq:joint convergence in distribution of ls estimate and other terms}.

\end{enumerate}

\subsection{Proof of Theorem \ref{thm:boundedness in probability of ls estimate}}

Based on Lemma \ref{lemma:convergence in distribution and tightness} and \ref{lemma:Relationships between modes of convergence}, we can derive Theorem \ref{thm:boundedness in probability of ls estimate} from Theorem \ref{thm:as and d convergence for terms of ls estimator} and \ref{thm:properties of ls estimate}.

\subsection{Proof of Theorems \ref{thm:boundedness of 4th moments of PVN}, \ref{thm:boundedness of 8th moment of PVN and 4th moment of VVN} and \ref{thm:boundedness of 8th moment of PPN and VVN}}

For $\E(\|{\Phi^{T}V}/{N}\|_{2}^4)$, there exists $\widetilde{M}_{1}>0$, irrespective of $N$, such that
\begin{align}\label{eq:derivation of boundedness of fourth moment of PVN}
\E\left(\left\|\frac{\Phi^{T}V}{N}\right\|_{2}^4\right)
=&\frac{1}{N^4}\E\left\{\sum_{i=1}^{n}\left[\sum_{t=1}^{N}u(t)v(t+i)\right]^2\right\}^2\nonumber\\
\leq &\frac{1}{N^4} 2^{n-1}\sum_{i=1}^{n}\left\{\E\left[\sum_{t=1}^{N}u(t)v(t+i)\right]^4 \right\}\nonumber\\
&[\text{using}\ \text{Lemma}\ \text{\ref{lem:C_r inequality}}]\nonumber\\
\leq & \frac{1}{N^2}\widetilde{M}_{1},
\end{align}
where the last step is derived from Assumptions \ref{asp:input}-\ref{asp:independence between input and noise} and Lemma \ref{lemma:strengthed boundedness of moments} with $m=4$. Similarly, under Assumption \ref{asp:input}-\ref{asp:bounded 8th moment of et and vt}, we can obtain \eqref{eq:boundedness of 8th moment of PVN} using Lemma \ref{lemma:strengthed boundedness of moments} with $m=8$. 

For $\E(\|\Phi^{T}\Phi/N-\Sigma\|_{F}^4)$, there exists $\widetilde{M}_{2}>0$, irrespective of $N$, such that
\begin{align}
&\E\left(\left\|\frac{\Phi^{T}\Phi}{N}-\Sigma\right\|_{F}^{4} \right)\nonumber\\
\leq & \E\Bigg\{\sum_{i=1}^{n}\sum_{j=1}^{n} \Bigg[\frac{1}{N}\sum_{t=1}^{N} u(t-i)u(t-j)\nonumber\\
&-\E\left(u(t-i)u(t-j)\right) \Bigg]^2\Bigg\}^2\nonumber\\
\leq& \frac{1}{N^4}2^{n^2-1}\sum_{i=1}^{n}\sum_{j=1}^{n}\E\left\{\sum_{t=1}^{N}\left[ u(t-i)u(t-j)\right.\right.\nonumber\\
&\left.-\E\left(u(t-i)u(t-j)\right) \right] \Bigg\}^4\nonumber\\
\leq& \frac{1}{N^2}\widetilde{M}_{2},
\end{align}
where the second step is derived using Lemma \ref{lem:C_r inequality} and the last step we apply Assumptions \ref{asp:input}, \ref{asp:bounded 8th moment of et and vt} and Lemma \ref{lemma:strengthed boundedness of moments} with $m=4$. Similarly, we can obtain \eqref{eq:boundedness of 8th moment of diff of PPN and Sigma} using Assumptions \ref{asp:input}, \ref{asp:bounded 16th moment of et and vt} and Lemma \ref{lemma:strengthed boundedness of moments} with $m=8$.

At the same time, there exists $\widetilde{M}_{4}>0$, irrespective of $N$, such that
\begin{align}\label{eq:middle result of boundedness of 4th moment of VVN diff}
\E\left|\frac{V^{T}V}{N}-\sigma^2\right|^4=& \frac{1}{N^4}\E\left\{ \sum_{t=1}^{N}\left[v^2(t)-\sigma^2\right]  \right\}^4\nonumber\\
\leq& \frac{1}{N^2}\widetilde{M}_{4},
\end{align}
where we apply Assumptions \ref{asp:noise}, \ref{asp:bounded 8th moment of et and vt} and Lemma \ref{lemma:strengthed boundedness of moments} with $m=4$. 
Similarly, we can obtain \eqref{eq:boundedness of 8th moment of diff of noise variance estimator} using Assumption \ref{asp:noise}, \ref{asp:bounded 16th moment of et and vt} and Lemma \ref{lemma:strengthed boundedness of moments} with $m=8$.

\subsection{Proof of Theorem \ref{thm:preliminary results about P and S}}

\begin{enumerate}
	
	\item For \eqref{eq:almost sure convergence of sqrtN_NPPinv}, we can apply \eqref{eq:almost sure convergence of PP_N} and Lemma \ref{lemma:continuous mapping theorem for vector case}.
	
	\item \underline{Proof of \eqref{eq:almost sure convergence of S_inv}, \eqref{eq:almost sure convergence of sqrtN_Sinv_Pinv} and \eqref{eq:difference of S_inv and P_inv}}
	
	For
	\begin{align}
	\hat{S}^{-1}=&[P+\widehat{\sigma^2}(\Phi^{T}\Phi)^{-1}]^{-1}.
	\end{align}
	using \eqref{eq:almost sure convergence of PP_N}, \eqref{eq:almost sure convergence of noise variance estimator} and Lemma \ref{lemma:continuous mapping theorem for vector case}, we can derive \eqref{eq:almost sure convergence of S_inv} as $N\to\infty$.
	The derivation of \eqref{eq:difference of S_inv and P_inv} is straightforward.
	Then for \eqref{eq:difference of S_inv and P_inv},
	we can apply \eqref{eq:almost sure convergence of PP_N}, \eqref{eq:almost sure convergence of noise variance estimator}, \eqref{eq:almost sure convergence of S_inv}, Lemma \ref{lemma:continuous mapping theorem for vector case} and \ref{lemma:Slutsky Theorem} to derive \eqref{eq:almost sure convergence of sqrtN_Sinv_Pinv}.
	
	\item \underline{Proof of \eqref{eq:1st derivatives of matrix inverse}, \eqref{eq:2nd derivatives of matrix inverse}, \eqref{eq:almost sure convergence of 1st order derivatives of S_inv and P_inv} and \eqref{eq:almost sure convergence of sqrtN difference of 1st order derivatives of S_inv and P_inv}}
	
	Using the basic identity of the inverse matrix derivative \cite[Page 9, Section 2.2, (59)]{PP2012}, we can derive \eqref{eq:1st derivatives of matrix inverse} and \eqref{eq:2nd derivatives of matrix inverse}.
	
	Then from \eqref{eq:1st derivatives of matrix inverse}, we have
	\begin{align}\label{eq:difference of 1st order derivatives of matrix inverse}
	&\frac{\partial \hat{S}^{-1}}{\partial \eta_{k}}-\frac{\partial P^{-1}}{\partial \eta_{k}}\nonumber\\
	=&(P^{-1}-\hat{S}^{-1})\frac{\partial P}{\partial \eta_{k}}P^{-1}
	+\hat{S}^{-1}\frac{\partial P}{\partial \eta_{k}}(P^{-1}-\hat{S}^{-1}),
	\end{align}
	which derives \eqref{eq:almost sure convergence of 1st order derivatives of S_inv and P_inv} as $N\to\infty$ based on Lemma \ref{lemma:Relationships between modes of convergence}, \ref{lemma:Slutsky Theorem} and \eqref{eq:almost sure convergence of S_inv}.
	
	Furthermore, combining \eqref{eq:difference of S_inv and P_inv} with \eqref{eq:difference of 1st order derivatives of matrix inverse}, we can derive \eqref{eq:almost sure convergence of sqrtN difference of 1st order derivatives of S_inv and P_inv} using \eqref{eq:almost sure convergence of PP_N}, \eqref{eq:almost sure convergence of noise variance estimator}, \eqref{eq:almost sure convergence of S_inv}, Lemma \ref{lemma:continuous mapping theorem for vector case} and \ref{lemma:Slutsky Theorem}.
	
	\item \underline{Proof of \eqref{eq:almost sure convergence of S_inv to P_inv at convergent estimate}}
	
	It can be known that there exists a constant $M>0$ such that $\|\hat{S}^{-1}\|_{F}<\|P^{-1}\|_{F}\leq M$ for any $\eta\in\widetilde{\Omega}$. Then we can obtain
	\begin{align}
	\sup_{\eta\in\widetilde{\Omega}}|\hat{S}^{-1}-P^{-1}|=&\sup_{\eta\in\widetilde{\Omega}}|\hat{S}^{-1}(\hat{S}-P)P^{-1}|\nonumber\\
	<&\widehat{\sigma^2} M^2 \frac{1}{N} \|N(\Phi^{T}\Phi)^{-1}\|_{F}\overset{a.s.}\to 0,
	\end{align}
	as $N\to\infty$, where the last step uses \eqref{eq:almost sure convergence of PP_N}, \eqref{eq:almost sure convergence of noise variance estimator}, Lemma \ref{lemma:continuous mapping theorem for vector case} and \ref{lemma:Slutsky Theorem}.
	It means that $\hat{S}^{-1}$ converges to $P^{-1}$ almost surely and uniformly in $\widetilde{\Omega}$.
	Further applying Lemma \ref{lemma:almost sure convergence of a function at convergent estimate} with $\hat{\eta}_{N}\overset{a.s.}\to \eta^{*}$ as $N\to\infty$, we have \eqref{eq:almost sure convergence of S_inv to P_inv at convergent estimate}.
	
	\item \underline{Proof of \eqref{eq:diff of S_inv and P_inv for general P} and \eqref{eq:equivalent form of diff of S_inv and P_inv}}
	
	According to \eqref{eq:def of hat_S_inv}, we have
	\begin{align}\label{eq:diff of Sinv at hat_etaN and Pinv at eta_star}
	&\hat{S}^{-1}(\hat{\eta}_{N})-P^{-1}(\eta^{*})\nonumber\\
	=&-\hat{S}^{-1}(\hat{\eta}_{N})\left[\hat{S}(\hat{\eta}_{N})-P(\eta^{*})\right]P^{-1}(\eta^{*})\nonumber\\
	=&-\hat{S}^{-1}(\hat{\eta}_{N})\left[P(\hat{\eta}_{N})-P(\eta^{*}) \right]P^{-1}(\eta^{*})\nonumber\\
	&-\widehat{\sigma^2}\hat{S}^{-1}(\hat{\eta}_{N})(\Phi^{T}\Phi)^{-1}P^{-1}(\eta^{*}).
	\end{align}
	Then we apply the first-order Taylor expansion of $P(\hat{\eta}_{N})$ at $\hat{\eta}_{N}=\eta^{*}$,
	\begin{align}\label{eq:Taylor expansion of P_hat_etaN}
	P(\hat{\eta}_{N})=&P(\eta^{*})
	+\sum_{k=1}^{p}\left.\frac{\partial P(\eta)}{\partial \eta_{k}}\right|_{\eta=\tilde{\eta}_{N}}(\hat{\eta}_{N,k}-\eta^{*}_{k})\nonumber\\
	=&P(\eta^{*})+\sum_{k=1}^{p}\left.\frac{\partial P(\eta)}{\partial \eta_{k}}\right|_{\eta=\tilde{\eta}_{N}}e_{k}^{T}(\hat{\eta}_{N}-\eta^{*}),
	\end{align}
	where we use the Lagrange's form of the remainder term, $\hat{\eta}_{N,k}$ and $\eta^{*}_{k}$ denote the $k$th element of $\hat{\eta}_{N}$ and $\eta^{*}$, respectively.
	Inserting \eqref{eq:Taylor expansion of P_hat_etaN} into \eqref{eq:diff of Sinv at hat_etaN and Pinv at eta_star}, we can obtain \eqref{eq:diff of S_inv and P_inv for general P}.
	
	In particular, for $P=\eta I_{n}$, the derivation of \eqref{eq:equivalent form of diff of S_inv and P_inv} is straightforward.
	
\end{enumerate}


\def\thesectiondis{\thesection.}                   
\def\thesubsectiondis{\thesection.\arabic{subsection}.}          
\def\thesubsubsectiondis{\thesection.\thesubsection.\arabic{subsubsection}.}

\renewcommand{\thesection}{B}
\renewcommand{\thesubsection}{B.\arabic{subsection}}
\setcounter{subsection}{0}

\setcounter{theorem}{0}
\setcounter{lemma}{0}
\setcounter{corollary}{0}
\renewcommand{\thelemma}{B.\arabic{lemma}}
\renewcommand{\thecorollary}{B.\arabic{corollary}}

\renewcommand{\theequation}{B.\arabic{equation}}
\setcounter{equation}{0}

\section*{Appendix B}\label{sec:Appendix B}

Fundamental lemmas and some preliminary results with their proofs are shown in Appendix B.

\subsection{Matrix Norm Inequalities }

\begin{lemma}\label{lemma:matrix norm inequality}(\cite{PP2012} Chapter 10.3 Page 61-62,\cite{HF2013} Page 68-72)
	 For $a\in\R^{m\times 1}$, $b\in\R^{m\times 1}$, $B\in\R^{m\times m}$ and $C\in\R^{m\times m}$, we have
		\begin{align}
		\label{eq:sum of 2 norm}
		\|a+b\|_{2}\leq &\|a\|_{2}+\|b\|_{2},\\
		\label{eq:inner product and 2 norm inequality}
		|a^{T}b|\leq & \|a\|_{2}\|b\|_{2},\\
		\label{eq:Frobenius norm and 2 norm}
		\|Bb\|_{2}\leq & \|B\|_{F}\|b\|_{2},\\
		\label{eq:Frobenius norm product}
		\|BC\|_{F}\leq&\|B\|_{F}\|C\|_{F}\\
		\label{eq:Frobenius norm sum}
		\|B+C\|_{F}\leq&\|B\|_{F}+\|C\|_{F}\\
		\label{eq:trace and frobenius norm inequality}
		|\Tr(B)|\leq&\sqrt{m}\|B\|_{F}.
		\end{align}
\end{lemma}

\pf Suppose that $a_{i}$ and $b_{i}$ are $i$th elements of $a$ and $b$, respectively. And let $b_{i,j}$ denote the $(i,j)$th element of $B$.
For \eqref{eq:inner product and 2 norm inequality}, it can be derived from
\begin{align}
|a^{T}b|=&\left|\sum_{i=1}^{m}a_{i}b_{i}\right|
\leq \sqrt{\sum_{i=1}^{m}a_{i}^2\sum_{j=1}^{m}b_{j}^2}=\|a\|_{2}\|b\|_{2},
\end{align}
which uses Cauchy-Schwarz inequaliy.

For \eqref{eq:Frobenius norm and 2 norm}, we have
\begin{align}
\|Bb\|_{2}
\leq& \|B\|_{2}\|b\|_{2}\ [\text{using}\ \text{properties}\ \text{of}\ \text{induced}\ \text{norm}]\nonumber\\
\leq& \|B\|_{F}\|b\|_{2},
\end{align}
where the last step comes from $\|B\|_{2}\leq\|B\|_{F}$ shown in Section $10.4.4$ in \cite{PP2012}.

For \eqref{eq:trace and frobenius norm inequality}, it can be seen that
\begin{align}
|\Tr(B)|=&\left|\sum_{i=1}^{m}B_{i,i}\right|\nonumber\\
\leq& \sqrt{m\sum_{i=1}^{m}|B_{i,i}|^2}\ \ [\text{using}\ \text{Cauchy-Schwarz}\ \text{Inequality}]\nonumber\\
\leq& \sqrt{m\sum_{i=1}^{m}\sum_{j=1}^{m}|B_{i,j}|^2}\nonumber\\
=&\sqrt{m}\|B\|_{F}.
\end{align}
%

\subsection{Ergodic Theory}

\begin{lemma}\label{lemma:ergodic theory}(\cite{Ljung1999} Page 43 Theorem 2.3)
 Let $\{s(t)\}$ be a stationary stochastic process with
		\begin{align}
		\E[s(t)]=&m(t)\\
		R_{s}(\tau)=&\E[s(t)s(t-\tau)].
		\end{align}
		Assume that
		\begin{align}
		s(t)-m(t)=x(t)=\sum_{k=0}^{\infty}h(k)e(t-k)=H(q)e(t),
		\end{align}
		where $\{e(t)\}$ is a sequence of independent random variables with zero mean values, $\E[e^{2}(t)]=\lambda_{t}$ and bounded fourth moments, and $H(q)$ is a stable filter, i.e. $\sum_{k=0}^{\infty}|h(k)|<\infty$. Then, as $N\to\infty$,
		\begin{align}
		&\frac{1}{N}\sum_{t=1}^{N}s(t)s(t-\tau)\overset{a.s.}\to  R_{s}(\tau)\\
		&\frac{1}{N}\sum_{t=1}^{N}\{s(t)m(t-\tau)-\E[s(t)m(t-\tau)]\} \overset{a.s.}\to 0.
		\end{align}
		This lemma is a special case of Theorem $2.3$ in \cite{Ljung1999}.
\end{lemma}

%

\subsection{Markov's Inequality}

\begin{lemma}\label{lem:Markov's inequality}(\cite{Gut2013} Page 120 Theorem 1.1)
 For a random variable $X\in\R$, suppose that $\E|X|^{r}<\infty$ for some $r>0$. Then for any $\epsilon>0$,
		\begin{align}
		\Pr(|X|>\epsilon)\leq \frac{\E|X|^{r}}{\epsilon^{r}}.
		\end{align}
		In particular, when $r=2$, it is also known as the Chebyshev's inequality.
\end{lemma}

\subsection{Borel-Cantelli's Lemma}

\begin{lemma}\label{lem:Borel-Cantelli's Lemma} (\cite{Ljung1999} Page 542 (\Rmnum{1}.18))
 For the random variable $X_{N}\in\R$, $\forall\epsilon>0$, if
		\begin{align}
		\sum_{N=1}^{\infty}\Pr(|X_{N}|>\epsilon)<\infty,
		\end{align}
		then as $N\to\infty$, we have
		\begin{align}
		X_{N}\overset{a.s.}\to 0.
		\end{align}
\end{lemma}

\subsection{Cauchy-Schwarz Inequality in Probability Theory}

\begin{lemma}\label{lem:cauchy schwarz inequality}(\cite{Gut2013} Page 130 Theorem 3.1)
	Suppose that $X_{1}\in\R$ and $X_{2}\in\R$ are random variables with finite variances. Then we have
		\begin{align}
		|\E(X_{1}X_{2})|\leq\E|X_{1}X_{2}|\leq\sqrt{\E(X_{1}^2)\E(X_{2}^2)}.
		\end{align}
\end{lemma}

\subsection{The $C_{r}$ Inequality}

\begin{lemma}\label{lem:C_r inequality}(\cite{Gut2013} Page 127 Theorem 2.2)
 For random variables $X_{1}\in\R$ and $X_{2}\in\R$, suppose that $\E|X_{1}|^{r}<\infty$ and $\E|X_{2}|^{r}<\infty$ for $r>0$. Then we have
		\begin{align}
		\E|X_{1}+X_{2}|^{r}\leq c_{r}(\E|X_{1}|^{r}+\E|X_{2}|^{r}),
		\end{align}
		where $c_{r}=1$ when $r\leq 1$ and $c_{r}=2^{r-1}$ when $r\geq 1$.
\end{lemma}

\subsection{Lyapounov Inequality}

\begin{lemma}\label{lemma:Lyapounov inequality}(\cite{Gut2013} Page 127 Theorem 2.5)
	For the random variable $X\in\R$ and $0<r\leq p$, we have
	\begin{align}
	\left(\E|X|^{r}\right)^{1/r}\leq \left(\E|X|^{p}\right)^{1/p}.
	\end{align}
\end{lemma}

\subsection{Center Limit Theorem (CLT) for $M$-dependent Sequence}

\begin{lemma}\label{lem:CLT for M-dependent sequence}(\cite{Ljung1999} Lemma 9.A1,\cite{Orey1958})
 Consider the sum of doubly indexed random vectors $\{x_{N}(t)\}_{t=1}^{N}$:
		\begin{align}
		Z_{N}=\sum_{t=1}^{N}x_{N}(t),
		\end{align}
		where $\E[x_{N}(t)]=0$. Suppose that $\{x_{N}(t)\}$ is $M-$dependent for an integer $M$, i.e. $\{x_{N}(1),x_{N}(2),\cdots,x_{N}(s)\}$ and $\{x_{N}(t),x_{N}(t+1),\cdots,x_{N}(n)\}$ are independent if $t-s>M$. In addition, assume that
		\begin{align}\label{eq:condition for 2nd moment}
		&\limsup\limits_{N\to\infty}\sum_{k=1}^{N}\E\|x_{N}(k)\|_{2}^2<\infty,\\
		\label{eq:condition for 2+delta moment}
		&\lim_{N\to\infty}\sum_{k=1}^{N}\E\|x_{N}(k)\|_{2}^{2+\delta}=0\ \text{for}\ \text{some}\ \delta>0.
		\end{align}
		Let
		\begin{align}
		Q_{z}=\lim_{N\to\infty}\E(Z_{N}Z_{N}^{T}).
		\end{align}
		Thus we have as $N\to\infty$,
		\begin{align}
		Z_{N}(t)\overset{d}\to \mathcal{N}(0,Q_{z}).
		\end{align}
\end{lemma}

\subsection{Limiting Theorem of Sum}

\begin{lemma}\label{lem:Limiting theorem of sum}(\cite{Ljung1999} Lemma 9.A2, \cite{DB1953})
 Let
		\begin{align}
		S_{N}=Z_{M}(N)+X_{M}(N),\ M,N=1,2,\cdots
		\end{align}
		such that
		\begin{subequations}
			\begin{align}
			\E\|X_{M}(N)\|_{2}^{2}\leq& C_{M},\ \lim_{M\to\infty}C_{M}=0\\
			\Pr\{Z_{M}(N)\leq z\}=&F_{M,N}(z)\\
			\lim_{N\to\infty}F_{M,N}(z)=&F_{M}(z)\\
			\lim_{M\to\infty}F_{M}(z)=&F(z),
			\end{align}
		\end{subequations}
		Then we have
		\begin{align}
		\lim_{N\to\infty}\Pr\{S_{N}\leq z\}=F(z).
		\end{align}
\end{lemma}

\subsection{Convergence in Distribution and Boundedness in Probability}

\begin{lemma}\label{lemma:convergence in distribution and tightness}(\cite{Olav2002} Lemma 3.8, \cite{Vaart1998} Page 8 Theorem 2.4)
	 Let $X_{N}\in\R$ be a random variable satisfying as $N\to\infty$
		\begin{align}
		X_{N}\overset{d}\to X.
		\end{align}
		Then $X_{N}$ is uniformly tight, which is also known to be bounded in probability, i.e. $X_{N}=O_{p}(1)$.	
\end{lemma}

%
%

\subsection{Continuous Mapping Theorem}

\begin{lemma}\label{lemma:continuous mapping theorem for vector case}(\cite{Vaart1998} Theorem $2.3$)
 Let $g:\R^{k}\mapsto\R^{m}$ be continuous at every point of a set $D$ such that $\text{Pr}(X\in D)=1$.
		\begin{enumerate}
			\item If $X_{N}\overset{d}\to X$ as $N\to\infty$, then $g(X_{N})\overset{d}\to g(X)$;
			\item If $X_{N}\overset{p}\to X$ as $N\to\infty$, then $g(X_{N})\overset{p}\to g(X)$;
			\item If $X_{N}\overset{a.s.}\to X$ as $N\to\infty$, then $g(X_{N})\overset{a.s.}\to g(X)$.
		\end{enumerate}
\end{lemma}

\subsection{Relationships between modes of convergence}

\begin{lemma}\label{lemma:Relationships between modes of convergence}(\cite{Vaart1998} Theorem $2.7$)
 Let $X_{N}$, $X$ and $Z_{N}$ be random vectors. Then as $N\to\infty$,
		\begin{enumerate}
			\item $X_{N}\overset{a.s.}\to X$ implies $X_{N}\overset{p}\to X$;
			\item $X_{N}\overset{p}\to X$ implies $X_{N}\overset{d}\to X$;
			\item $X_{N}\overset{p}\to c$ for a constant $c$ if and only if $X_{N}\overset{d}\to c$.
		\end{enumerate}
\end{lemma}

\subsection{Slutsky Theorem}

\begin{lemma}\label{lemma:Slutsky Theorem}(\cite{Vaart1998} Theorem $2.8$)
 Let $X_{N}$, $X$ and $Y_{N}$ be random vectors or variables. If as $N\to\infty$, $X_{N}\overset{d}\to X$ and $Y_{N}\overset{d}\to c$ for a constant $c$, then
		\begin{enumerate}
			\item $X_{N}+Y_{N}\overset{d}\to X+c$;
			\item $X_{N}Y_{N}\overset{d}\to cX$;
			\item \label{item:equivalence of convergence in distribution and probability to a constant} $Y_{N}^{-1}X_{N}\overset{d}\to c^{-1}X$ provided $c\not=0$.	
		\end{enumerate}
\end{lemma}

\subsection{The Kolmogorov Strong Law}

\begin{lemma}\label{lem:kolmogorov strong law}\cite[Page 295, Theorem 6.1(a)]{Gut2013}
 For the random sequence $\{X_{i}\}_{i=1}^{N}$ with $X_{i}\in\R$ for $i=1,2,\cdots,N$, if $\E|X_{i}|<\infty$ and $\E(X_{i})=\mu$, as $N\to\infty$, we have
		\begin{align}
		\frac{1}{N}\sum_{i=1}^{N}X_{i}\overset{a.s.}\to \mu.
		\end{align}
\end{lemma}

\subsection{Squeeze Theorem for Convergence in Distribution}

\begin{lemma}\label{lemma:squeeze theorem}
	{\it
		Let $X_{N}\in\R$, $Y_{N,1}\in\R$, $Y_{N,2}\in\R$ and $X\in\R$ be random. If as $N\to\infty$, $Y_{N,1}\overset{d.}\to X$, $Y_{N,2}\overset{d.}\to X$, and $X_{N}$ satisfies $Y_{N,1}\leq X_{N}\leq Y_{N,2}$, we have
		\begin{align}
		X_{N}\overset{d.}\to X,
		\end{align}
		as $N\to\infty$.
	}
\end{lemma}

\pf According to the definition of the convergence in distribution, we can know that
\begin{align}
\lim_{N\to\infty}\Pr(Y_{N,1}\leq x)=&\Pr(X\leq x)\\
\lim_{N\to\infty}\Pr(Y_{N,2}\leq x)=&\Pr(X\leq x),
\end{align}
for every $x$ at which the distribution function $F_{X}:x\mapsto\Pr(X\leq x)$ is continuous. Meanwhile, $Y_{N,1}\leq X_{N}\leq Y_{N,2}$ derives
\begin{align}
\Pr(Y_{N,1}\leq x)\leq \Pr(X_{N}\leq x) \leq \Pr(Y_{N,2}\leq x).
\end{align}
Then based on the standard Squeeze theorem (\cite[Page 104 Theorem 3.3.6]{Hou2014}), we have $X_{N}\overset{d.}\to X$ as $N\to\infty$.
\qed

\subsection{Inequalities of Determinant and Trace}

\begin{lemma}\label{lem:relationship between log determinant and trace}
	{\it
		For a positive definite matrix $A\in\R^{m\times m}$, we have
		\begin{align}\label{eq:inequality of log determinant and trace}
		\Tr(I_{n}-A^{-1})\leq \log\det(A) \leq \Tr(A-I_{n}).
		\end{align}
	}
\end{lemma}

\pf Our first step is to show that for any $x>0$,
\begin{align}
f_{1}(x)=&1-\frac{1}{x}-\log(x)\leq 0\\
f_{2}(x)=&\log(x)-x+1\leq 0.
\end{align}
The first order derivatives of $f_{1}(x)$ and $f_{2}(x)$ with respect to $x$, respectively are
\begin{align}
\frac{d f_{1}(x)}{d x}=&\frac{1}{x^2}(1-x)=0\ \Rightarrow\ x=1,\\
\frac{d f_{2}(x)}{d x}=&\frac{1}{x}(1-x)=0\ \Rightarrow\ x=1,
\end{align}
which derives that $\forall\ x>0$, $f_{1}(x)\leq f_{1}(1)=0$ and $f_{2}(x)\leq f_{2}(1)=0$.

Then we rewrite \eqref{eq:inequality of log determinant and trace} using eigenvalues of $A$:
\begin{align*}
\Tr(I_{n}-A^{-1})-\log\det(A)=&\sum_{i=1}^{m}\left[1-\frac{1}{\lambda_{i}(A)}-\log(\lambda_{i}(A))\right],\\
\log\det(A)-\Tr(A-I_{n})=&\sum_{i=1}^{m}\left[\log(\lambda_{i}(A))-\lambda_{i}(A)+1\right].
\end{align*}
Since $A$ is positive definite, we have $\lambda_{i}(A)>0$. Then applying $f_{1}(\lambda_{i}(A))\leq 0$ and $f_{2}(\lambda_{i}(A))\leq 0$ can lead to \eqref{eq:inequality of log determinant and trace}. \qed

\subsection{Almost Sure Convergence of Convergent Function at Convergent Estimate}

\begin{lemma}\label{lemma:almost sure convergence of a function at convergent estimate}
	{\it Suppose that as $N\to\infty$, $M_{N}(\eta)$ converges almost surely to a non-stochastic function $M(\eta)$ uniformly in a compact set $D$ containing $\eta^{*}$ and $\hat{\eta}$. If $\hat{\eta} \overset{a.s.}\to \eta^{*}$ as $N\to\infty$, and $M(\eta)$ is continuous in $D$, we have as $N\to\infty$,
		\begin{align}
		M_{N}(\hat{\eta})\overset{a.s.}\to M(\eta^{*}).
		\end{align}
	}	
\end{lemma}

\pf	It can be known that
\begin{align*}
|M_{N}(\hat{\eta})-M(\eta^{*})|=&|M_{N}(\hat{\eta})-M(\hat{\eta})+M(\hat{\eta})-M(\eta^{*})|\nonumber\\
\leq& |M_{N}(\hat{\eta})-M(\hat{\eta})|+|M(\hat{\eta})-M(\eta^{*})|,
\end{align*}
where we need to consider two terms $|M_{N}(\hat{\eta})-M(\hat{\eta})|$ and $|M(\hat{\eta})-M(\eta^{*})|$, respectively.

For the first term, since $M_{N}(\eta)$ converges almost surely to a non-stochastic function $M_{\eta}$ uniformly in a compact set $D$, we have
\begin{align}
|M_{N}(\hat{\eta})-M(\hat{\eta})|
\leq \sup_{\eta\in D}|M_{N}(\eta)-M(\eta)|\overset{a.s.}\to 0,
\end{align}
which comes from the almost surely uniform convergence of $M_{N}(\eta)$ to $M_{\eta}$ in $D$.

For the second term, since $M(\eta)$ is continuous and $\hat{\eta} \overset{a.s.}\to \eta^{*}$, we can apply Lemma \ref{lemma:continuous mapping theorem for vector case} to obtain as $N\to\infty$,
\begin{align}
|M(\hat{\eta})-M(\eta^{*})|\overset{a.s.}\to 0.
\end{align}

Consequently, as $N\to\infty$, $M_{N}(\hat{\eta})\overset{a.s.}\to M(\eta^{*})$.
\qed

\subsection{Moments of Filtered White Noise}

\begin{lemma}\cite[(A.48)-(A.52)]{shumway2000}\label{lemma:properties of covariance function of filtered white noise}
 For a stationary process with $t=1,2,\cdots$,
		\begin{align}
		x(t)=\sum_{k=0}^{\infty}h(k)w(t-k),
		\end{align}
		where $w(t)$ is independent and identically distributed (i.i.d.) with zero mean, variance $\sigma_{w}^2$ and fourth moment
		\begin{align}
		\E[w^4(t)]=\epsilon\sigma_{w}^4<\infty,\ \epsilon\in\R\ \text{is}\ \text{some}\ \text{constant},
		\end{align}
		and $\sum_{k=0}^{\infty}|h(k)|<\infty$, for fixed $\tau,\tau'\in\Z$,
		if we define the covariance function
		\begin{align}
		\gamma(\tau)=&\E[x(t)x(t+\tau)]\ \text{with}\ \gamma(-\tau)=\gamma(\tau),
		\end{align}
		and its sample estimator 
		\begin{align}
		\widehat{\gamma}(\tau)=&\frac{1}{N}\sum_{t=1}^{N}x(t)x(t+\tau),
		\end{align}
		then we have
		\begin{align}
		&\gamma(\tau)=\sigma_{w}^2\sum_{k=0}^{\infty}h(k)h(k+|\tau|),\nonumber\\
		&\lim_{N\to\infty}N\left[\widehat{\gamma}(\tau)-\gamma(\tau)\right]\left[\widehat{\gamma}(\tau')-\gamma(\tau')\right]\nonumber\\
		=&(\epsilon-3)\gamma(\tau)\gamma(\tau')+\sum_{i=-\infty}^{\infty}\left[\gamma(i)\gamma(i+\tau-\tau')+\gamma(i+\tau)\gamma(i-\tau')\right].
		\end{align}
\end{lemma}

\subsection{Relationship of Convergence in Distribution and Convergence of Moments}

\begin{lemma}\label{lem:relationship of convergence in dist and convergence of moments for scalars}(\cite[Page 100 Theorem 4.5.2]{Chung2001})
		For random $x_{N}\in\R$ and $x\in\R$, if $x_{N}\overset{d}\to x$ as $N\to\infty$, and for some $p>0$, $\exists\ M>0$ such that
		\begin{align}
		\sup_{N}\E\left(|x_{N}|^{p}\right)\leq M,
		\end{align}
		then for each $0<r<p$ and $r\in\N$, we have
		\begin{align}
		\lim_{N\to\infty}\E\left(x_{N}^{r}\right)=\E\left(x^{r}\right).
		\end{align}
\end{lemma}

Moreover, Lemma \ref{lem:relationship of convergence in dist and convergence of moments for scalars} can be expended to the multidimensional case to obtain the following results.

\begin{lemma}\label{lem:relationship of convergence in dist and convergence of expectation for multidim}
	For random $X_{N}\in\R^{m\times 1}$ and $X\in\R^{m\times 1}$, if $X_{N}\overset{d}\to X$ as $N\to\infty$, and $\exists\ M>0$, where $M$ is irrespective of $N$, such that
	\begin{align}
	\sup_{N}\E\|X_{N}\|_{2}^{2}\leq M_{1},
	\end{align}
	then we have
	\begin{align}\label{eq:convergence of 1st moment}
	\lim_{N\to\infty}\E[X_{N}]=\E[X].
	\end{align}
\end{lemma}

\pf Define that the $i$th elements of $X_{N}$ and $X$ with $i=1,\cdots,m$ are $X_{N,i}$ and $X_{i}$, respectively. Since $\exists\ M_{1}>0$ such that
\begin{align}
\sup_{N}\E\|X_{N}\|_{2}^2=\sup_{N}\E\left( \sum_{i=1}^{m}X_{N,i}^2 \right)\leq M_{1},
\end{align}
we can see that for each $i=1,\cdots,m$,
\begin{align}
\sup_{N}\E\left( X_{N,i}^2 \right)\leq \sup_{N}\E\|X_{N}\|_{2}^2\leq M_{1}.
\end{align}
At the same time, since as $N\to\infty$, $X_{N}\overset{d}\to X$, we have as $N\to\infty$, for $i=1,\cdots,m$,
\begin{align}
X_{N,i}\overset{d}\to X_{i}.
\end{align}
Then applying Lemma \ref{lem:relationship of convergence in dist and convergence of moments for scalars} with
\begin{align}
x_{N}=X_{N,i},\ x=X_{i},\ p=2,\ r=1,\ i=1,\cdots,m,
\end{align}
we can obtain \eqref{eq:convergence of 1st moment}.
\qed

\begin{lemma}\label{lem:relationship of convergence in dist and convergence of moments for multidim}
	 For random $X_{N}\in\R^{m\times 1}$ and $X\in\R^{m\times 1}$, if $X_{N}\overset{d}\to X$ as $N\to\infty$, and $\exists\ M>0$, where $M$ is irrespective of $N$, such that
			\begin{align}
			\sup_{N}\E\|X_{N}\|_{2}^{4}\leq M,
			\end{align}
			then we have
			\begin{align}\label{eq:convergence of covariance matrix}
			\lim_{N\to\infty}\E(X_{N}X_{N}^{T})=\E(XX^{T}).
			\end{align}
\end{lemma}

\pf Define that the $i$th elements of $X_{N}$ and $X$ with $i=1,\cdots,m$ are $X_{N,i}$ and $X_{i}$, respectively. It can be seen that the derivation of \eqref{eq:convergence of covariance matrix} is equivalent to showing that for every $i,j=1,\cdots,m$,
\begin{align}
\label{eq:middle result convergence of covariance matrix element}
\lim_{N\to\infty} \E(X_{N,i}X_{N,j})=&\E(X_{i}X_{j}).
\end{align}
Since
\begin{align}
\sup_{N}\E\|X_{N}\|_{2}^{4}=\sup_{N}\E\left[\left(\sum_{i=1}^{m}X_{N,i}^2\right)^{2}\right]
\end{align}
is bounded, it can be derived that for each $i=1,\cdots,m$,
\begin{align}\label{eq:middle boundedness result 1}
\sup_{N}\E\left(X_{N,i}^4\right)\leq& \sup_{N}\E\|X_{N}\|_{2}^{4}\leq M_{2},
\end{align}
which leads to that for each $i,j=1,\cdots,m$,
\begin{align}
\sup_{N}\E(|X_{N,i}X_{N,j}|^{2})
\leq &\sup_{N}\sqrt{\E(X_{N,i}^{4})\E(X_{N,j}^{4})}\nonumber\\
&[\text{using}\ \text{Lemma}\ \text{\ref{lem:cauchy schwarz inequality}}]\nonumber\\
\label{eq:middle boundedness result 2}
\leq&  M_{2}.
\end{align}
Then applying \eqref{eq:middle boundedness result 2} and Lemma \ref{lem:relationship of convergence in dist and convergence of moments for scalars} with
\begin{align}
x_{N}=X_{N,i}X_{N,j},\ x=X_{i}X_{j},\ p=2,\ r=1,\ i,j=1,\cdots,m,
\end{align}
we can obtain \eqref{eq:middle result convergence of covariance matrix element}, which leads to \eqref{eq:convergence of covariance matrix}. 
\qed

\subsection{Relationship of Essentially Boundedness and Boundedness of Moments}

\begin{lemma}\label{lem:relationship of essentially boundedness and bounded expectation for positive rv}
		For the random variable $X_{N}>0$, if $X_{N}$ is essentially bounded, i.e. there exits $M>0$ such that for all $N=1,2,\cdots$, the subset of event space $\Omega_{1}=\{\omega||X_{N}|>M\}$ satisfies
		\begin{align}
			\Pr(\Omega_{1})=&0,
			\Pr(\overline{\Omega_{1}})=1,
		\end{align}
		where $\overline{\Omega_{1}}$ is the complement of $\Omega_{1}$, then we have
		\begin{align}
		\E(X_{N}^{k})\leq M^{k},\ \forall\ k\geq 1,\ k\in\N.
		\end{align}
\end{lemma}

\pf We consider that $X_{N}$ is a continuous random variable (the discrete case is similar). It follows that
\begin{align}
\E(X_{N}^{k})=&\int_{\Omega_{1}}X_{N}^{k}dP(\omega)+\int_{\overline{\Omega_{1}}}X_{N}^{k}dP(\omega)\nonumber\\
=&\int_{\overline{\Omega_{1}}}X_{N}^{k}dP(\omega)\nonumber\\
\leq& \int_{\overline{\Omega_{1}}}M^{k}dP(\omega)\nonumber\\
=&M^{k}\Pr(\overline{\Omega_{1}})\nonumber\\
=&M^{k}.
\end{align}
\qed

\subsection{Upper and Lower Bounds of A Trace}

\begin{lemma}\label{lemma:bounds of traces}
	For $A\in\R^{m_{1}\times m_{2}}$, $B\in\R^{m_{1}\times m_{1}}$ and $k\in\Z^{+}$, if $B$ is positive definite, define that the largest and smallest eigenvalue of $B$ are $\lambda_{1}(B)$ and $\lambda_{m_{1}}(B)$, respectively. Let $u_{B,m_{1}}\in\R^{m_{1}}$ denote the eigenvector associated with $\lambda_{m_{1}}(B)$ and $\cond(B)$ denote the condition number of $B$, defined as $\cond(B)=\lambda_{1}(B)/\lambda_{m_{1}}(B)$.
	\begin{itemize}
		\item If $A$ is irrespective of $B$ and $u_{B,m_{1}}^{T}A\not=0$, then there exists $B_{L},B_{U}>0$, irrespective of $\cond(B)$, such that
		\begin{align}\label{eq:bounds of trace1}
		\frac{B_{L}}{\lambda_{1}^{k}(B)}{\cond}^{k}(B)\leq \Tr(A^{T}&B^{-k}A) \leq \frac{B_{U}}{\lambda_{1}^k(B)}{\cond}^{k}(B),
		\end{align}
		where
		\begin{align}\label{eq:explicit expressions of bounds of trace1}
		B_{L}=u_{B,m_{1}}^{T}AA^{T}u_{B,m_{1}},\
		B_{U}=\Tr(AA^{T}).
		\end{align}
		
		\item If $m_{1}=m_{2}$, $A$ is irrespective of $B$, and $u_{B,m_{1}}^{T}Au_{B,m_{1}}\not=0$, there exists $B_{L},B_{U}>0$, irrespective of $\cond(B)$, such that
		\begin{align}
		\label{eq:bounds of trace2}
		\frac{B_{L}}{\lambda_{1}^{3}(B)}{\cond}^{3}(B)\leq& \Tr(B^{-1}A^{T}B^{-1}AB^{-1})\nonumber\\
		\leq& \frac{B_{U}}{\lambda_{1}^3(B)}{\cond}^{3}(B),
		\end{align}
		where 
		\begin{align}\label{eq:explicit expressions of bounds of trace2}
		B_{L}=\left(u_{B,m_{1}}^{T}Au_{B,m_{1}}\right)^2,\
		B_{U}=\Tr(AA^{T}).
		\end{align}
	\end{itemize}
\end{lemma}

\pf Define the eigenvalue decomposition (EVD) of $B$ as follows,
\begin{align}
\label{eq:EVD of B}
B=&U_{B}\Lambda_{B}U_{B}^{T}=\sum_{i=1}^{m_{1}}\lambda_{i}(B)u_{B,i}u_{B,i}^{T},
\end{align}
where $\Lambda_{B}\in\R^{m_{1}\times m_{1}}$ is diagonal with $\lambda_{1}(B)\geq\cdots\geq\lambda_{m_{1}}(B)>0$, eigenvalues of $B$, and $U_{B}\in\R^{m_{1}\times m_{1}}$ is an orthogonal matrix with $u_{B,i}\in\R^{m_{1}}$ being its $i$th column vector. Let $\cond(B)=\lambda_{1}(B)/\lambda_{m_{1}}(B)$ denote the condition number of $B$.

For \eqref{eq:bounds of trace1}, inserting \eqref{eq:EVD of B}, we have
\begin{align}
&\Tr(A^{T}B^{-k}A)\nonumber\\
=&\Tr\left[A^{T}\sum_{i=1}^{m_{1}}\frac{1}{\lambda_{i}^k(B)}u_{B,i}u_{B,i}^{T}A\right]\nonumber\\
=&\frac{1}{\lambda_{1}^k(B)}\sum_{i=1}^{m_{1}}\frac{\lambda_{1}^k(B)}{\lambda_{i}^k(B)}u_{B,i}^{T}AA^{T}u_{B,i}\nonumber\\
=&\frac{1}{\lambda_{1}^k(B)}\left[u_{B,1}^{T}AA^{T}u_{B,1} +\sum_{i=2}^{m_{1}}\frac{\lambda_{1}^k(B)}{\lambda_{i}^k(B)}u_{B,i}^{T}AA^{T}u_{B,i} \right],
\end{align}
Then if we let
\begin{align}
B_{L}=&u_{B,m_{1}}^{T}AA^{T}u_{B,m_{1}},\\
B_{U}=&\sum_{i=1}^{m_{1}}u_{B,i}^{T}AA^{T}u_{B,i}\nonumber\\
=&\Tr\left(AA^{T}\sum_{i=1}^{m_{1}}u_{B,i}u_{B,i}^{T}\right)\nonumber\\
=&\Tr(AA^{T}),
\end{align}
where
\begin{align}
\sum_{i=1}^{m_{1}}u_{B,i}u_{B,i}^{T}=U_{B}U_{B}^{T}=I_{m_{1}},
\end{align}
we can obtain \eqref{eq:bounds of trace1}. 

For $m_{1}=m_{2}$, the derivation of \eqref{eq:bounds of trace2} is analogous. Inserting \eqref{eq:EVD of B}, we have
\begin{align}
&\Tr(B^{-1}A^{T}B^{-1}AB^{-1})\nonumber\\
=&\Tr\left[\sum_{i=1}^{m_{1}}\sum_{j=1}^{m_{1}}\frac{1}{\lambda_{i}(B)\lambda_{j}^2(B)}A^{T}u_{B,i}u_{B,i}^{T}Au_{B,j}u_{B,j}^{T}\right]\nonumber\\
=&\frac{1}{\lambda_{1}^3(B)}\sum_{i=1}^{m_{1}}\sum_{j=1}^{m_{1}}\frac{\lambda_{1}^3(B)}{\lambda_{i}(B)\lambda_{j}^2(B)}
\left(u_{B,i}^{T}Au_{B,j}\right)^2.
\end{align}
Then if we let
\begin{align}
B_{L}=&\left(u_{B,m_{1}}^{T}Au_{B,m_{1}}\right)^2,\\
B_{U}=&\sum_{i=1}^{m_{1}}\sum_{j=1}^{m_{1}}\left(u_{B,i}^{T}Au_{B,j}\right)^2\nonumber\\
=&\sum_{i=1}^{m_{1}}u_{B,i}^{T}A\left(\sum_{j=1}^{m_{1}}u_{B,j}u_{B,j}^{T}\right)A^{T}u_{B,i}\nonumber\\
=&\Tr\left(AA^{T}\sum_{i=1}^{m_{1}}u_{B,i}u_{B,i}^{T}\right)\nonumber\\
=&\Tr(AA^{T}),
\end{align}
we can obtain \eqref{eq:bounds of trace2}

\qed

\subsection{Relationship between Noise Variance and SNR}

\begin{lemma}\label{lemma:relationship between noise variance and snr for fir model}
	For the FIR model \eqref{eq:vector-matrix form linear regression model}, as $N\to\infty$, we have
	\begin{align}\label{eq:as convergence of snr and noise variance}
	\text{SNR}\overset{a.s.}\to \frac{\theta_{0}^{T}\Sigma\theta_{0}}{\sigma^2}.
	\end{align}
\end{lemma}

\pf According to the definition of SNR, for the fixed noise variance, we mainly consider the almost sure convergence of the sample variance of noise-free outputs. For $i=1,\cdots,N$, let
\begin{align}
X_{i}=&\phi(i)^{T}\theta_{0}=\sum_{j=1}^{n}u(i-j)g_{j}^{0}.
\end{align}
The sample variance of noise-free outputs can be represented as
\begin{align}
&\frac{1}{N-1}\sum_{i=1}^{N}\left( X_{i}-\frac{1}{N}\sum_{j=1}^{N}X_{j} \right)^2\nonumber\\
=&\frac{N}{N-1}\left( \frac{1}{N}\sum_{i=1}^{N}X_{i}^2 \right)-\frac{N}{N-1}\left( \frac{1}{N}\sum_{i=1}^{N}X_{i} \right)^2.
\end{align}
Since as $N\to\infty$,
\begin{align}
\frac{1}{N}\sum_{i=1}^{N}X_{i}=&\sum_{j=1}^{n}g_{j}^{0}\frac{1}{N}\sum_{i=1}^{N}u(i-j)\overset{a.s.}\to 0,\\
\frac{1}{N}\sum_{i=1}^{N}X_{i}^2=&\sum_{j=1}^{n}\sum_{k=1}^{n}g_{j}^{0}g_{k}^{0}\frac{1}{N}\sum_{i=1}^{N}u(i-j)u(i-k)\nonumber\\
\overset{a.s.}\to & \sum_{j=1}^{n}\sum_{k=1}^{n}g_{j}^{0}g_{k}^{0}R_{u}(|j-k|)=\theta_{0}^{T}\Sigma\theta_{0},
\end{align}
which can be derived using Lemma \ref{lem:kolmogorov strong law}, then we can obtain \eqref{eq:as convergence of snr and noise variance}. \qed

\subsection{Some Preliminary Results}

\begin{lemma}\label{lem:Boundedness of Expectation}
 Let
		\begin{align}\label{eq:sequence w1}
		w_{1}(t)=&\sum_{k=0}^{\infty}\alpha_{1}(k)x(t-k),
		\end{align}
		where $\sum_{k=0}^{\infty}|\alpha_{1}(k)|<\infty$ and $\{x(t)\}$ is independent with zero mean.
		\begin{enumerate}
			\item If $\{x(t)\}$ is equipped with bounded fourth moments, we have
			\begin{align}
			\label{eq:boundedness 1 of norm expectation}
			&\E\left\{\sum_{t=1}^{N}\left[w_{1}(t-\tau_{1})w_{1}(t-\tau_{2})\right.\right.\nonumber\\
			&-\left.\E(w_{1}(t-\tau_{1})w_{1}(t-\tau_{2}))\right]\bigg\}^{2}
			\leq 4C_{\alpha}^4C_{x4}N,
			\end{align}
			where
			\begin{align}
			C_{\alpha}=&\sum_{k=0}^{\infty}|\alpha_{1}(k)|,\
			C_{x4}=\sup_{t}\E[x^{4}(t)].
			\end{align}
			
			\item If $\{x(t)\}$ is equipped with bounded second moments and there exists a random sequence $\{w_{2}(t)\}$ independent of $\{x(t)\}$, which is also independent with zero mean and bounded second moments, we have
			\begin{align}
			\label{eq:boundedness 2 of norm expectation}
			&\E\left[\sum_{t=1}^{N}w_{1}(t)w_{2}(t)\right]^{2}\leq C_{\alpha}^2C_{x2}C_{w}N,
			\end{align}
			where
			\begin{align}
			C_{x2}=\sup_{t}\E[x^{2}(t)],\
			C_{w}=\sup_{t}\E[w_{2}^{2}(t)].
			\end{align}
		\end{enumerate}	
\end{lemma}

\pf	The inequality \eqref{eq:boundedness 1 of norm expectation} can be derived from Lemma 2B.1 in \cite{Ljung1999}.

For \eqref{eq:boundedness 2 of norm expectation}, it can be known that
\begin{align}
&\E\left[\sum_{t=1}^{N}w_{1}(t)w_{2}(t)\right]^{2}\nonumber\\
=&\E\Bigg[\sum_{t=1}^{N}\sum_{l=1}^{N}\sum_{k_{1}=0}^{\infty}\sum_{k_{2}=0}^{\infty}\alpha_{1}(k_{1})x(t-k_{1})\nonumber\\
&w_{2}(t)\alpha_{1}(k_{2})x(l-k_{2})w_{2}(l)\Bigg]\nonumber\\
=&\sum_{t=1}^{N}\sum_{k_{1}=0}^{\infty}\alpha_{1}^2(k_{1})\E[x^2(t-k_{1})]\E[w_{2}^{2}(t)]\nonumber\\
\leq&\sum_{t=1}^{N}\left[\sum_{k_{1}=0}^{\infty}|\alpha_{1}(k_{1})|\right]^2 \sup_{t}\E[x^2(t-k_{1})] \sup_{t}\E[w_{2}^{2}(t)]\nonumber\\
\leq &C_{\alpha}^2C_{x2}C_{w}N.
\end{align}
where the last third step applies the respective independence of $\{x(t)\}$ and $\{w_{2}(t)\}$ and their mutual independence.
\qed

The following lemma is a strengthened version of Lemma \ref{lem:Boundedness of Expectation}.

\begin{lemma}\label{lemma:strengthed boundedness of moments}
	For $m=4$ or $m=8$, let
	\begin{align}
	w_{1}(t)=&\sum_{k=0}^{\infty}\alpha_{1}(k)x(t-k),
	\end{align}
	where $\sum_{k=0}^{\infty}|\alpha_{1}(k)|<\infty$ and $\{x(t)\}$ is independent with zero mean.
	\begin{enumerate}
		\item If $\{x(t)\}$ is equipped with bounded $2m$ moments, there exists $\widetilde{C}_{1}>0$, which is irrespective of $N$, such that
				\begin{align}
				\label{eq:boundedness of m moment of w1}
				&\E\left\{\sum_{t=1}^{N}\left[w_{1}(t-\tau_{1})w_{1}(t-\tau_{2})\right.\right.\nonumber\\
				&-\left.\E(w_{1}(t-\tau_{1})w_{1}(t-\tau_{2}))\right]\bigg\}^{m}
				\leq \widetilde{C}_{1}N^{m/2}.
				\end{align}
		\item If $\{x(t)\}$ is equipped with bounded $m$ moments and there exists a random sequence $\{w_{2}(t)\}$ independent of $\{x(t)\}$, which is also independent with zero mean and bounded $m$ moments, there exists $\widetilde{C}_{2}>0$, which is irrespective of $N$, such that
				\begin{align}
				\label{eq:boundedness of 2m moment of w1 and w2}
				&\E\left[\sum_{t=1}^{N}w_{1}(t)w_{2}(t)\right]^{m}\leq \widetilde{C}_{2}N^{m/2}.
				\end{align}
	\end{enumerate}
\end{lemma}

\pf Here we only derive the case of $m=4$, and the case of $m=8$ is similar and thus omitted. The main thought is analogous to the proof of Lemma \ref{lem:Boundedness of Expectation}.

For \eqref{eq:boundedness of m moment of w1} with $m=4$, without loss of generality (W.L.O.G.), we take $\tau_{1}=\tau_{2}=0$. Define that
\begin{align}
S_{w_{1}}^{N}=&\sum_{t=1}^{N}\left[w_{1}(t)^2-\E(w_{1}^2(t))\right]\\
=&\sum_{t=1}^{N}\sum_{k=0}^{\infty}\sum_{l=0}^{\infty}\alpha_{1}(k)\alpha_{1}(l)\beta(t,k,l),
\end{align}
where
\begin{align}
\beta(t,k,l)=\left[x(t-k)x(t-l)-\E[x^2(t-k)]\delta_{k,l}\right],
\end{align}
and $\delta_{k,l}$ denotes the Kronecker delta, which is equal to one if $k=l$ otherwise zero. Then 
\begin{align}
\left(S_{w_{1}}^{N}\right)^4
=&\sum_{t_{1}=1}^{N}\cdots\sum_{t_{4}=1}^{N}\sum_{k_{1}=0}^{\infty}\sum_{l_{1}=0}^{\infty}\cdots\sum_{k_{4}=0}^{\infty}
\sum_{l_{4}=0}^{\infty}\nonumber\\
&\alpha_{1}(k_{1})\alpha_{1}(l_{1})\cdots\alpha_{1}(k_{4})\alpha_{1}(l_{4})\nonumber\\
&\gamma(t_{1},\cdots,t_{4},k_{1},l_{1},\cdots,k_{4},l_{4}),
\end{align}
where
\begin{align}
\gamma(t_{1},\cdots,t_{4},k_{1},l_{1},\cdots,k_{4},l_{4})=\beta(t_{1},k_{1},l_{1})\cdots\beta(t_{4},k_{4},l_{4}).
\end{align}
Since $\{x(t)\}$ is independent, we can see that $\E(\gamma)$ is equal to zero unless at least some indices in $\beta(t_{i},k_{i},l_{i})$ with $i=1,\cdots,4$ coincide.
Thus, we have
\begin{align}
&\E\left\{\sum_{t=1}^{N}\left[w_{1}(t-\tau_{1})w_{1}(t-\tau_{2})\right.\right.\nonumber\\
&-\left.\E(w_{1}(t-\tau_{1})w_{1}(t-\tau_{2}))\right]\bigg\}^{m}\nonumber\\
\leq & \sum_{k_{1}=0}^{\infty}|\alpha_{1}(k_{1})|\sum_{l_{1}=0}^{\infty}|\alpha_{1}(l_{1})|
\sum_{k_{2}=0}^{\infty}|\alpha_{1}(k_{2})|\sum_{l_{2}=0}^{\infty}|\alpha_{1}(l_{2})|\nonumber\\
&\sum_{t_{1}=1}^{N}\sum_{t_{2}=1}^{N}16\sup_{t}\E\left[x^8(t)\right]\nonumber\\
\leq& \widetilde{C}_{1}N^2,
\end{align}
where
\begin{align}
\widetilde{C}_{1}\triangleq 16 \left[\sum_{k=0}^{\infty}|\alpha_{1}(k)| \right]^4 \sup_{t}\E\left[x^8(t)\right].
\end{align}

For \eqref{eq:boundedness of 2m moment of w1 and w2}, we have
\begin{align}
&\E\left[\sum_{t=1}^{N}w_{1}(t)w_{2}(t)\right]^{4}\nonumber\\
=&\sum_{t_{1}=1}^{N}\sum_{t_{2}=1}^{N}\sum_{k_{1}=0}^{\infty}\sum_{k_{2}=0}^{\infty}\alpha_{1}^2(k_{1})\alpha_{1}^2(k_{2})\nonumber\\
&\E\left[x^{2}(t_{1}-k_{1}) x^2(t_{2}-k_{2})\right]\E\left[w_{2}^2(t_{1}) w_{2}^2(t_{2})\right]\nonumber\\
\leq&N^2\left[\sum_{k_{1}=0}^{\infty}|\alpha_{1}(k_{1})|\right]^{4}\sup_{t}\E[x^{4}(t)]\sup_{t}\E[w_{2}^4(t)]\nonumber\\
\leq& \widetilde{C}_{2}N^2,
\end{align}
where the first step is based on the independence of $\{x(t)\}$ and $\{w_{2}(t)\}\}$ and their mutual independence, and
\begin{align}
\widetilde{C}_{2}\triangleq
\left[\sum_{k_{1}=0}^{\infty}|\alpha_{1}(k_{1})|\right]^{4}\sup_{t}\E[x^{4}(t)]\sup_{t}\E[w_{2}^4(t)].
\end{align}
\qed

\begin{lemma}\label{lemma:Almost Sure Convergence for Statistics with Bounded Expectation}
	{\it If there exists a constant $C_{s}$ such that
		\begin{align}\label{eq:upper bounds of expectation of sum of st}
		\E\left[\sum_{t=1}^{N}\tilde{s}(t)\right]^2\leq C_{s}N,
		\end{align}
		as $N\to\infty$, we have
		\begin{align}\label{eq:almost sure convergence of sample average}
		\frac{1}{N}\sum_{t=1}^{N}\tilde{s}(t)\overset{a.s.}\to 0.
		\end{align}
	}
\end{lemma}

\pf	Similar proof is included in the proof of Theorem 2B.2 in \cite{Ljung1999}.
Our proof consists of two steps.
\begin{itemize}
	\item Firstly we show that $\frac{1}{N^2}\sum_{t=1}^{N^2}\tilde{s}(t)\overset{a.s.}\to 0$ as $N\to\infty$.
	
	Since $\E[\sum_{t=1}^{N}\tilde{s}(t)]^2\leq C_{s}\cdot N$, it gives that
	\begin{align}
	\E\left[\frac{1}{N^2}\sum_{t=1}^{N^2}\tilde{s}(t)\right]^2\leq \frac{CN^2}{N^4}=\frac{C_{s}}{N^2}.
	\end{align}
	Based on the Chebyshev's inequality as shown in Lemma \ref{lem:Markov's inequality} with $r=2$, for any $\epsilon>0$, we have
	\begin{align}
	\text{Pr}\left(\left|\frac{1}{N^2}\sum_{t=1}^{N^2}\tilde{s}(t)\right|>\epsilon \right)
	\leq \frac{1}{\epsilon^2}\E\left[\frac{1}{N^2}\sum_{t=1}^{N^2}\tilde{s}(t)\right]^2\leq \frac{C_{s}}{\epsilon^2 N^2},
	\end{align}
	which yields that
	\begin{align}
	\sum_{N=1}^{\infty}\text{Pr}\left(\left|\frac{1}{N^2}\sum_{t=1}^{N^2}\tilde{s}(t)\right|\right)\leq \frac{C_{s}}{\epsilon}\sum_{N=1}^{\infty}\frac{1}{N^2}<\infty.
	\end{align}
	Using Lemma \ref{lem:Borel-Cantelli's Lemma} (Borel-Cantelli's lemma), it implies that as $N\to\infty$,
	\begin{align}\label{eq:almost sure convergence of N_squared sum of s_tau}
	\frac{1}{N^2}\sum_{t=1}^{N^2}\tilde{s}(t)\overset{a.s.}\to 0.
	\end{align}
	
	\item Secondly, we show that $\sup_{N^2\leq k\leq (N+1)^2}\frac{1}{k}\left|\sum_{t=1}^{k}\tilde{s}(t)\right|\overset{a.s.}\to 0$ as $N\to\infty$, which proves Lemma \ref{lemma:Almost Sure Convergence for Statistics with Bounded Expectation}.
	
	Suppose that $\sup_{N^2\leq k\leq (N+1)^2}\frac{1}{k}\left|\sum_{t=1}^{k}\tilde{s}(t)\right|$ is achieved at $k=k_{N}$. Then we can obtain that
	\begin{align}
	&\sup_{N^2\leq k\leq (N+1)^2}\frac{1}{k}\left|\sum_{t=1}^{k}\tilde{s}(t)\right|\nonumber\\
	=&\frac{1}{k_{N}}\left|\sum_{t=1}^{k_{N}}\tilde{s}(t)\right|\nonumber\\
	\label{eq:two terms of sup sum of s_tau}
	\leq&\frac{1}{k_{N}}\left|\sum_{t=1}^{N^2}\tilde{s}(t)\right|+\frac{1}{k_{N}}\left|\sum_{t=N^2+1}^{k_{N}}\tilde{s}(t)\right|.
	\end{align}
	\begin{itemize}
		\item[$-$] For the first term of \eqref{eq:two terms of sup sum of s_tau}, since $k_{N}\geq N^2$, it is known that
		\begin{align}
		\frac{1}{k_{N}}\left|\sum_{t=1}^{N^2}\tilde{s}(t)\right|\leq \frac{1}{N^2}\left|\sum_{t=1}^{N^2}\tilde{s}(t)\right|.
		\end{align}
		Combining with \eqref{eq:almost sure convergence of N_squared sum of s_tau}, we know that as $N\to\infty$,
		\begin{align}\label{eq:1st term of sup sum of s_tau}
		\frac{1}{k_{N}}\left|\sum_{t=1}^{N^2}\tilde{s}(t)\right|\overset{a.s.}\to 0.
		\end{align}
		
		\item[$-$] For the second term of \eqref{eq:two terms of sup sum of s_tau}, we have
		\begin{align}
		&\E\left[\frac{1}{k_{N}}\sum_{t=N^2+1}^{k_{N}}\tilde{s}(t)\right]^2\nonumber\\
		&\leq\frac{1}{N^4}\E\left[\sum_{t=N^2+1}^{k_{N}}\tilde{s}(t)\right]^2\ [\text{since}\ k_{N}\geq N^2]\nonumber\\
		&\leq\frac{1}{N^4}C_{s}[(N+1)^2-N^2-1](k_{N}-N^2-1)\nonumber\\ 
		&[\text{since}\ \eqref{eq:upper bounds of expectation of sum of st}]\nonumber\\
		&\leq\frac{1}{N^4}C_{s}[(N+1)^2-N^2-1]^2\ [\text{since}\ k_{N}\leq (N+1)^2]\nonumber\\
		&=\frac{4C}{N^2}.
		\end{align}
		Similarly, according to the Chebyshev's inequality (Lemma \ref{lem:Markov's inequality} with $r=2$) and Borel-Cantelli's lemma (Lemma \ref{lem:Borel-Cantelli's Lemma}), we can derive that as $N\to\infty$,
		\begin{align}\label{eq:2nd term of sup sum of s_tau}
		\frac{1}{k_{N}}\left|\sum_{t=N^2+1}^{k_{N}}\tilde{s}(t)\right|\overset{a.s.}\to 0.
		\end{align}
	\end{itemize}
	Consequently, inserting \eqref{eq:1st term of sup sum of s_tau} and \eqref{eq:2nd term of sup sum of s_tau} into \eqref{eq:two terms of sup sum of s_tau}, we show that as $N\to\infty$,
	\begin{align}
	\sup_{N^2\leq k\leq (N+1)^2}\frac{1}{k}\left|\sum_{t=1}^{k}\tilde{s}(t)\right|\overset{a.s.}\to 0,
	\end{align}
	which implies \eqref{eq:almost sure convergence of sample average}.
\end{itemize}
\qed

Then we can derive the following result using Lemma \ref{lem:cauchy schwarz inequality} and \ref{lem:C_r inequality}.
\begin{lemma}\label{lemma:bounded result}
	{\it Let $\beta_{i}\in\R$ and $\gamma_{i}\in\R$ with $i=1,\cdots,m$ be random variables and have bounded moments with order $4+\delta$ for some $\delta>0$. For
		\begin{align}
		X=\left[\begin{array}{c}\beta_{1}\gamma_{1}-\E(\beta_{1}\gamma_{1})\\ \vdots \\ \beta_{m}\gamma_{m}-\E(\beta_{m}\gamma_{m}) \end{array}\right],
		\end{align}
		there exists a constant $C_{\beta\gamma}$ satisfying
		\begin{align}
		\E\|X\|_{2}^{2+\delta}\leq C_{\beta\gamma}.
		\end{align}
		In particular, if $\beta_{i}$ and $\gamma_{i}$ are mutually independent, as long as they have bounded moments of order $2+\delta$ for some $\delta>0$, $\E\|X\|_{2}^{2+\delta}$ is still bounded.
	}
\end{lemma}

\pf	\begin{align}
&\E\|X\|_{2}^{2+\delta}\nonumber\\
=&\E\left\{\sum_{i=1}^m[\beta_{i}\gamma_{i}-\E(\beta_{i}\gamma_{i})]^2\right\}^{(2+\delta)/2}\nonumber\\
\leq&2^{\delta m/2}\sum_{i=1}^{m}\E|\beta_{i}\gamma_{i}-\E(\beta_{i}\gamma_{i})|^{2+\delta}\nonumber\\
&[\text{using}\ \text{Lemma}\ \ref{lem:C_r inequality}\ \text{with}\ r=(2+\delta)/2]\nonumber\\
\label{eq:middle result of upper bounds of 2_delta moments}
\leq&2^{\delta m/2+1+\delta}\sum_{i=1}^{m}[\E|\beta_{i}\gamma_{i}|^{2+\delta}+|\E(\beta_{i}\gamma_{i})|^{2+\delta}]\nonumber\\
&[\text{using}\ \text{Lemma}\ \text{\ref{lem:C_r inequality}}\ \text{with}\ r=2+\delta].
\end{align}
Apply the Cauchy-Schwarz inequality as shown in Lemma \ref{lem:cauchy schwarz inequality} to obtain
\begin{subequations}\label{eq:upper bounds of two terms about beta and gamma}
	\begin{align}\label{eq:upper boundes of term 1}
	\E|\beta_{i}\gamma_{i}|^{2+\delta}\leq& [\E(\beta_{i})^{4+2\delta}\E(\gamma_{i})^{4+2\delta}]^{1/2},\\
	|\E(\beta_{i}\gamma_{i})|^{2+\delta}\leq&[\E(\beta_{i})^2\E(\gamma_{i})^2]^{(2+\delta)/2}.
	\end{align}
\end{subequations}
Note that for a random variable $X$, if $\E|X|^{s}<\infty$, for any $0<r\leq s$, we have $\E|X|^{r}<\infty$. Combining \eqref{eq:middle result of upper bounds of 2_delta moments} and \eqref{eq:upper bounds of two terms about beta and gamma}, since $\beta_{i}$ and $\gamma_{i}$ have bounded moments of order $4+\delta$ (or lower order) and \eqref{eq:middle result of upper bounds of 2_delta moments} is a finite sum of bounded moments, hence $\E\|X\|_{2}^{2+\delta}$ is also bounded for some fixed $\delta$.

Define that
\begin{align}
C_{\beta}=&\sup_{i}\E(\beta_{i})^2,\ \tilde{C}_{\beta}=\sup_{i}\E(\beta_{i})^{4+\delta},\\
C_{\gamma}=&\sup_{i}\E(\gamma_{i})^2,\ \tilde{C}_{\gamma}=\sup_{i}\E(\gamma_{i})^{4+\delta}.
\end{align}
Thus we have
\begin{align}
\E\|X\|_{2}^{2+\delta}\leq C_{\beta\gamma},
\end{align}
where
\begin{align}
C_{\beta\gamma}=2^{\delta m/2+1+\delta}m[(\tilde{C}_{\beta}\tilde{C}_{\gamma})^{1/2}+(C_{\beta}C_{\gamma})^{(2+\delta)/2}].
\end{align}

In particular, if $\beta_{i}$ and $\gamma_{i}$ are mutually independent, \eqref{eq:upper boundes of term 1} can be adjusted as
\begin{align}
\E|\beta_{i}\gamma_{i}|^{2+\delta}=\E|\beta_{i}|^{2+\delta}\E|\gamma_{i}|^{2+\delta},
\end{align}
which will still be bounded if $\beta_{i}$ and $\gamma_{i}$ have bounded moments of order $2+\delta$ for some $\delta>0$.
\qed


\bibliographystyle{unsrt}
\bibliography{database}

\end{document}